\documentclass[reqno]{amsart}
\usepackage{amssymb,amsthm,fullpage}
\usepackage{enumerate}
\usepackage{cite}
\usepackage{soul}
\usepackage{amsmath}
\usepackage{amsfonts}
\usepackage{fancyhdr,xcolor}
\usepackage{tikz}
\usetikzlibrary{calc}
\usepackage{wrapfig}
\usepackage{listings}
\newtheorem{Thm}{Theorem}[section]
\newtheorem{Lem}[Thm]{Lemma}
\newtheorem{rem}[Thm]{Remark}
\newtheorem{prop}[Thm]{Proposition}
\newtheorem{Cor}[Thm]{Corollary}

\newcommand{\mc}{\mathcal}
\newcommand{\mb}{\mathbb}

\newcommand{\zetabar}{\overline{\zeta}}
\newcommand{\Z}{\mathbb{Z}}
\newcommand{\R}{\mathbb{R}}
\newcommand{\mF}{\mathcal{F}}
\newcommand{\nablab}{\nabla_b}

\DeclareMathOperator{\Tr}{Tr}
\newcommand{\eps}{\epsilon}

\newcommand{\norma}[1]{{\left\vert\kern-0.25ex\left\vert\kern-0.25ex\left\vert #1 
    \right\vert\kern-0.25ex\right\vert\kern-0.25ex\right\vert}}
\newcommand{\norm}[1]{{\left\vert\kern-0.25ex\left\vert #1 
    \right\vert\kern-0.25ex\right\vert}}
    
\DeclareMathOperator{\e}{\mathrm{err}}

\DeclareMathOperator{\Real}{Re}
\DeclareMathOperator{\Imag}{Im}

\newcommand{\defn}[1]{\emph{#1}}

\def\sideremark#1{\ifvmode\leavevmode\fi\vadjust{\vbox to0pt{\vss
 \hbox to 0pt{\hskip\hsize\hskip1em
 \vbox{\hsize3cm\tiny\raggedright\pretolerance10000
 \noindent #1\hfill}\hss}\vbox to8pt{\vfil}\vss}}}

\numberwithin{equation}{section}

\begin{document}
\title[The Neumann problem on the Clifford torus in $\mathbb{S}^3$]{The Neumann problem on the domain in $\mathbb{S}^3$ bounded by the Clifford torus}
\author{Jeffrey S.\ Case}
\address{Department of Mathematics \\ Penn State University \\ University Park, PA 16802, USA}
\email{jscase@psu.edu}
\author{Eric Chen}
\address{Department of Mathematics \\ University of California \\ Berkeley, CA 94720-3840, USA}
\email{ecc@berkeley.edu}
\author{Yi Wang}
\address{Mathematics Department \\ Johns Hopkins University \\ Baltimore, MD, 21218, USA}
\email{ywang261@jhu.edu}
\author{Paul Yang}
\address{Department of Mathematics \\ Princeton University \\ Princeton, NJ 08544-1000, USA}
\email{yang@math.princeton.edu}
\author{Po-Lam Yung}
\address{ Australian National University, Canberra, ACT 2601, Australia
\& The Chinese University of Hong Kong, Shatin, Hong Kong}
\email{polam.yung@anu.edu.au, \, plyung@math.cuhk.edu.hk}

\keywords{subelliptic PDE; Neumann problem; CR Yamabe problem; layer potentials}
\subjclass[2020]{Primary 58J32; Secondary 35R03, 35H20}
\begin{abstract}
We discuss the solution of the Neumann problem associated with the CR Yamabe operator on a subset $\Omega$ of the CR manifold $\mathbb{S}^3$ bounded by the Clifford torus $\Sigma$.
We also discuss the Yamabe-type problem of finding a contact form on $\Omega$ which has zero Tanaka--Webster scalar curvature and for which $\Sigma$ has constant $p$-mean curvature.
\end{abstract}
\dedicatory{Dedicated to David Jerison on the occasion of his 70th birthday}
\maketitle

\section{Introduction}\label{sec_intro}


We view $\mathbb{S}^3=\{(\zeta^1,\zeta^2)\in\mathbb{C}^2:~|\zeta^1|^2+|\zeta^2|^2 = 1\}$ as the boundary of the unit ball in $\mathbb{C}^2$, and equip it with a pseudohermitian structure associated to the contact form
\[
\theta:= i \overline{\partial} (|\zeta|^2-1) = \frac{i}{2} (\zeta^1 d\zetabar^1 + \zeta^2 d\zetabar^2 - \zetabar^1 d\zeta^1 - \zetabar^2 d\zeta^2).
\] 
Let $\Omega\subset \mathbb{S}^3$ be the domain
\begin{equation}\label{eq:Omega_def}
\Omega=\{(\zeta^1,\zeta^2) \in \mathbb{S}^3 \colon |\zeta^1|<|\zeta^2|\},
\end{equation}
whose boundary is given by
\[
\Sigma:=\partial\Omega=\{(\zeta^1,\zeta^2) \in \mathbb{S}^3 \colon |\zeta^1|=|\zeta^2|=\frac{1}{\sqrt{2}}\},
\]
the Clifford torus. 

On $\Omega$, we give an explicit solution to the Neumann problem, \begin{align}
\begin{cases}
\label{N}
L u=0\quad &\text{in }\Omega,
\\
\nabla_\nu u=h\quad&\text{on }\partial\Omega .
\end{cases}
\end{align}
Here $L$ denotes the CR Yamabe operator~\eqref{eqn:cr-yamabe-operator-sphere} on $(\mathbb{S}^3,\theta)$ and $\nabla_\nu$ denotes the one-sided horizontal normal derivative for a function $u \in C^{\infty}(\overline{\Omega})$,
i.e.
\[
\nabla_{\nu} u(\zeta) := \lim_{t \to 0^+} \frac{u(\gamma_{\zeta}(t)) - u(\zeta)}{t}, \quad \zeta \in \partial \Omega,
\]
where $\gamma_{\zeta}(t)$ is a curve so that $\gamma_{\zeta}(0) =\zeta$ and $\gamma'(0)$ is the inward horizontal normal $\nu$ to $\partial \Omega$ at $\zeta \in \partial \Omega$.

Initially we assume that our boundary data $h$ is smooth on $\Sigma$, i.e.\ $h \in C^{\infty}(\partial \Omega)$, and solve \eqref{N}.
Our solution is given  in terms of a suitable single layer potential $\mathcal{S}$, defined for $f \in C^1(\Sigma)$ by
\[
\mathcal{S}f(\zeta) := \oint_{\eta \in \Sigma} G(\zeta,\eta) f(\eta) d\sigma(\eta), \quad \zeta \in \overline{\Omega} ,
\]
where $G(\zeta,\eta)$ is the Green's function of the CR Yamabe operator $L$, and $d\sigma$ is the surface measure on $\Sigma$ corresponding to the volume form $\iota_{\nu} \theta \wedge d\theta$.

\begin{Thm} \label{thm:Neumann}
If $h \in C^{\infty}(\Sigma)$, then there exists a unique $u \in C^{\infty}(\overline{\Omega})$ that solves the Neumann problem \eqref{N}.
Furthermore, $u$ is given by $\mathcal{S}f$ on $\overline{\Omega}$ for some $f \in C^{\infty}(\Sigma)$, where $f$ is determined by $h$ from
\begin{equation} \label{eq:f_to_h}
f = (-\frac{1}{2}I + \mathcal{N})^{-1} h,
\end{equation}
and $\mathcal{N}$ is the singular integral operator given by
\[
\mathcal{N}f(\zeta) := p.v. \oint_{\eta \in \Sigma} (\nabla_{\nu})_{\zeta} G(\zeta,\eta) f(\eta) d\sigma(\eta), \quad \zeta \in \Sigma.
\]
\end{Thm}
Some of the arguments follow those in \cite{OV}, where \eqref{N} is solved on flag domains of $\mathbb{H}^1$.
In contrast to the Euclidean case~\cite{MR501367}, the operator $\mathcal{N}$ is not a compact operator on $L^2(\Sigma)$; it is only bounded on $L^2(\Sigma)$, but not smoothing of any positive order. As a result, a more careful analysis is necessary to show that $-\frac{1}{2}I + \mathcal{N}$ is invertible on $L^2(\Sigma)$ (and all higher order $L^2$-based Sobolev spaces).

A general result of Nhieu~\cite{Nhieu2001} implies the existence of a solution $u$ to~\eqref{N} in the Folland--Stein~\cite{FollandStein1974} space $S^{1,2}(\Omega)$ of functions $u \in L^2(\Omega)$ whose horizontal gradient is also in $L^2(\Omega)$.
We recover this result using the single layer potential.

\begin{Thm} \label{thm:NeumannL2}
If $h \in L^2(\Sigma)$, then there exists a unique $u \in S^{1,2}(\Omega)$ such that \eqref{N} is satisfied in the weak sense, i.e.
\[
\int_{\Omega} \left(\langle \nabla_b u , \nabla_b \phi \rangle + \frac{R}{4} \, u \, \phi \right) \, \theta \wedge d\theta + \oint_{\Sigma} h \phi \, d\sigma = 0
\]
for all $\phi \in C^{\infty}(\overline{\Omega})$.
Furthermore, $u$ is given by $\mathcal{S}f$ on $\Omega$ where $f \in L^2(\Sigma)$ is determined from $h$ by \eqref{eq:f_to_h}.
\end{Thm}

The key observation in the proof of Theorem~\ref{thm:NeumannL2} is that the single layer potential extends to a bounded linear map $\mathcal{S} \colon L^2(\Sigma) \to S^{1,2}(\Omega)$.

Our work is motivated in part by the desire to formulate and study the CR boundary Yamabe problem (cf. \cite{Escobar1992a}):
Given a closed CR three-manifold with boundary $(M^3,T^{1,0})$, construct a Webster-flat contact form with respect to which the boundary has constant $p$-mean curvature $H$~\cite{ChengHwangMalchiodiYang2005}.
Given a contact form $\theta$, the contact form $\widehat{\theta} := u^2\theta$ satisfies these properties if and only if $u$ is a positive critical point of the functional $\mc{F} \colon \mc{V} \to \mb{R}$,
\begin{align*}
 \mc{F}(u) & := \int_M \left( \lvert\nabla_b u\rvert^2 + \frac{R}{4}u^2 \right) \, \theta \wedge d\theta + \frac{1}{3}\oint_{\partial M} Hu^2\, d\sigma , \\
 \mc{V} & := \left\{ u \in C^\infty(M) \mathrel{}:\mathrel{} \oint_{\partial M} \lvert u\rvert^3 \, d\sigma = 1 \right\} .
\end{align*}
One way to construct such a contact form is to show that there is a smooth, positive function which realizes the CR boundary Yamabe constant
\begin{equation*}
 Y(M,T^{1,0}) := \inf \left\{ \mc{F}(u) \mathrel{}:\mathrel{} u \in \mc{V} \right\} .
\end{equation*}

If $\partial M$ has no characteristic points, then the Sobolev trace embedding theorem~\cite{Nhieu2001} implies that the restriction map $C^\infty(\overline{M}) \ni u \mapsto u\rvert_{\partial M} \in C^\infty(\partial M)$ extends to a continuous linear map $\Tr \colon S^{1,2}(M) \to L^2(\partial M)$; in particular, $\mathcal{F}$ is well-defined on $S^{1,2}(M)$.
As in the Riemannian case~\cite{Escobar1992corr}, it holds that $Y(M,T^{1,0})>-\infty$ if and only if the Dirichlet eigenvalues of $L$ are positive;
i.e.
\begin{equation*}
 \lambda_{1,D}(L) := \inf \left\{ \mc{F}(u) \mathrel{}:\mathrel{} u\rvert_{\partial M} = 0, \int_M \lvert u\rvert^2 \, \theta \wedge d\theta = 1 \right\} > 0 .
\end{equation*}
The positivity of $\lambda_{1,D}(L)$ also implies that the first Steklov eigenvalue
\begin{equation}
 \label{eqn:first-neumann}
 \mu_{1}(L) := \inf \left\{ \mc{F}(u) \mathrel{}:\mathrel{} \oint_{\partial M} \lvert u \rvert^2 \, d\sigma = 1 \right\}
\end{equation}
of $L$ is finite. 

In Section~\ref{sec:compute-boundary-form} we discuss the equivalence of the signs of $\mu_{1}(L)$ and $Y(M,T^{1,0})$, under the assumption that minimizers of $\mu_{1}(L)$ are smooth up to the boundary. Such regularity assumption can be verified in the case where $M = \overline{\Omega}$, where $\Omega$ is the domain in $\mathbb{S}^3$ defined by \eqref{eq:Omega_def}. Note that $H = 0$ for the Clifford torus $\Sigma = \partial \Omega$.

\begin{prop}
 \label{prop:clifford-steklov-regularity}
 Let $(\Omega,\theta)$ be the interior of the Clifford torus with the standard spherical contact form.
 If $u \in S^{1,2}(\Omega)$ minimizes $\mu_{1}(L)$, then $u \in C^\infty(\overline{\Omega})$.
\end{prop}

Proposition~\ref{prop:clifford-steklov-regularity} will be proved by observing that a minimizer $u \in S^{1,2}(\Omega)$ of $\mu_1(L)$ is a weak solution of
\[
\begin{cases}
Lu = 0 , & \text{on $\Omega$} \\
\nabla_{\nu} u = \mu \Tr(u) , &\text{in $\partial \Omega$} ,
\end{cases}
\]
where $\mu = -\mathcal{F}(u)$ is constant, and then using tools developed for the proof of Theorem~\ref{thm:NeumannL2}.




Since the standard contact form on $\mathbb{S}^3$ has positive Tanaka--Webster scalar curvature and is such that the Clifford torus is $p$-minimal, the Sobolev trace embedding theorem~\cite{Nhieu2001} implies that the boundary Yamabe constant $Y(\Omega,T^{1,0})$ is positive.
We identify an explicit critical point of the functional $\mc{F} \colon \mc{V} \to \mathbb{R}$.

\begin{Thm}
 \label{thm:clifford-extremals}
 Let $(\Omega,T^{1,0},\theta)$ be the interior of the Clifford torus.
 Set
 \begin{equation*}
  u := \frac{{}_2F{}_1\bigl( \frac{1}{2} , \frac{1}{2} ; 1 ; \lvert z^2\rvert^2 \bigr)}{{}_2F{}_1\bigl( \frac{1}{2} , \frac{1}{2} ; 1 ; \frac{1}{2} \bigr)},
 \end{equation*}
 where
 \begin{align*}
  {}_2F{}_1(a,b;c;x) & := \sum_{n=0}^\infty \frac{(a)_n(b)_n}{(c)_nn!}x^n , \\
  (a)_0 & := 1 , \\
  (a)_n & := a(a+1)\dotsm(a+n-1) , && \text{if $n \geq 1$},
 \end{align*}
 is the standard Gaussian hypergeometric function.
 Then $u^2\theta$ is a scalar flat contact form on $\Omega$ with respect to which $\Sigma$ has constant $p$-mean curvature.
\end{Thm}

The solution of Theorem~\ref{thm:clifford-extremals} is normalized so that $\theta$ and $u^2\theta$ coincide on $\Sigma$.


The article is organized as follows. In Section~\ref{sect2}, we give some background in pseudohermitian geometry and Fourier analysis, in the specific setting of the Clifford torus. In Section~\ref{sect3}, we compute the single layer potential in local coordinates. Theorems~\ref{thm:Neumann}, \ref{thm:NeumannL2}, Proposition~\ref{prop:clifford-steklov-regularity} and Theorem~\ref{thm:clifford-extremals} will be proved in Sections~\ref{sect4}, \ref{sect5}, \ref{sect6} and \ref{sec:compute-boundary-form}, respectively.

\subsection*{Acknowledgements}
Case is partially supported by the Simons Foundation (Grant \#524601).
Chen is partially supported by NSF Award DMS-3103392.
Wang is partially supported by NSF Career award DMS-1845033.
Yang is partially supported by the Simons Foundation (Grant \#1006518).
Yung is partially supported by a Future Fellowship FT200100399 from the Australian Research Council.

\section{Setup and notations} \label{sect2}

\subsection{Pseudohermitian geometry} 
\label{subsec:ph-geometry}
A \defn{CR three-manifold} is a pair $(M^3,T^{1,0})$ consisting of a real three-manifold $M$ and a complex rank $1$ distribution $T^{1,0} \subset T_{\mathbb{C}}M$.
Let $H := \mathrm{Re}\bigl( T^{1,0} \oplus \overline{T^{0,1}}\bigr)$ denote the space of \emph{horizontal vectors}.
Then
\begin{equation*}
 J(Z + \overline{Z}) := iZ - i\overline{Z} ,
\end{equation*}
$Z \in T^{1,0}$, defines an integrable almost complex structure on $H$.
We say that $(M^3,T^{1,0})$ is \emph{nondegenerate} if locally there is a real one-form $\theta$ such that $\ker\theta = H$ and $\theta \wedge d\theta$ is nowhere-vanishing.
Nondegenerate CR three-manifolds are orientable~\cite[Lemma~19]{Jacobowitz1990}, and hence there is a global real one-form $\theta$ such that $\ker\theta = H$;
in this case we call $\theta$ a \defn{contact form}.

A \defn{(strictly pseudoconvex) pseudohermitian manifold} is a triple $(M^3,T^{1,0},\theta)$ consisting of a nondegenerate CR three-manifold $(M^3,T^{1,0})$ and a contact form $\theta$ such that $d\theta(Z,\overline{Z})>0$ for all nonzero $Z \in T^{1,0}$.
The \defn{Reeb vector field} is the unique vector field $T$ such that $\theta(T) = 1$a and $d\theta(T,\cdot) = 0$.

An \defn{admissible coframe} for $(M^3,T^{1,0},\theta)$ is a nowhere-vanishing local complex-valued one-form $\theta^1$ such that $\theta^1(T)=0$ and $\theta^1(\overline{Z})=0$ for all $Z \in T^{1,0}$.
Set $\theta^{\bar1} := \overline{\theta^1}$.
Then $\{ \theta, \theta^1 , \theta^{\bar1} \}$ is a local coframe for $T_{\mathbb{C}}^\ast M$.
It follows that
\begin{equation*}
 d\theta = ih_{1\bar1} \, \theta^1 \wedge \theta^{\bar1} .
\end{equation*}
Note that $h_{1\bar1}>0$.
Let $\{ T , Z_1 , Z_{\bar1} \}$ be the dual frame to $\{ \theta, \theta^1 , \theta^{\bar1} \}$.
This (globally) determines a positive definite inner product on $H$ by
\begin{equation*}
 \langle \Real a^1Z_1 , \Real b^1Z_1 \rangle :=  \frac{1}{2}\Real h_{1\bar 1}a^1\overline{b^1} .
\end{equation*}


Let $\theta^1$ be an admissible coframe for $(M^3,T^{1,0},\theta)$.
Then there is a unique complex-valued one-form $\omega_1{}^1$ such that
\begin{equation*}
 \begin{cases}
  d\theta^1 = \theta^1 \wedge \omega_1{}^1 + A^1{}_{\bar1} \, \theta \wedge \theta^{\bar 1} , \\
  dh_{1\bar1} = \omega_1{}^1h_{1\bar1} + \omega_{\bar1}{}^{\bar1}h_{1\bar1} ,
 \end{cases}
\end{equation*}
where $\omega_{\bar1}{}^{\bar1} := \overline{\omega_1{}^1}$.
The \defn{Tanaka--Webster connection} $\nabla$ is uniquely determined from $\nabla Z_1 := \omega_1{}^1 \otimes Z_1$ and $\nabla T := 0$ by complex linearity.
The \defn{pseudohermitian torsion} is the globally-defined tensor $A_{11}\,\theta^1\otimes\theta^1$.
The \defn{Tanaka--Webster scalar curvature} is the globally-defined function $R$ determined by
\begin{equation*}
 d\omega_1{}^1 = Rh_{1\bar1} \, \theta^1 \wedge \theta^{\bar 1} \mod \theta .
\end{equation*}
We say that $\theta$ is \defn{scalar flat} if $R=0$.

Given a function $f \in C^\infty(M)$, we denote by $\nablab f$ the \emph{subgradient} of $f$; i.e.\ the restriction of $df$ to $H$.
The \defn{sublaplacian} $\Delta_b \colon C^\infty(M) \to C^\infty(M)$ is
\begin{equation*}
 \Delta_b := \nablab^\ast\nablab ,
\end{equation*}
where $\nablab^\ast$ is the formal $L^2$-adjoint of $\nablab$ with respect to $\theta \wedge d\theta$.
Locally,
\begin{equation*}
 \Delta_b u = -h^{\bar{1}1} \left( (Z_1 Z_{\bar{1}} + Z_{\bar{1}} Z_1) u - \omega_1{}^1(Z_{\bar{1}}) Z_1 u - \omega_{\bar{1}}{}^{\bar{1}} (Z_1) Z_{\bar{1}} u \right) .
\end{equation*}
It is readily computed (cf.\ \cite[Equation~(2.4)]{Lee1986}) that
\begin{equation}
 \label{eqn:ibp1}
 \int_M \langle \nablab u , \nablab w \rangle \, \theta \wedge d\theta = \int_M u\Delta_bw \, \theta \wedge d\theta + 2\Real \oint_{\partial M} iu(Z_1w) \, \theta \wedge \theta^1 .
\end{equation}

Suppose that $(M^3,T^{1,0},\theta)$ is a pseudohermitian three-manifold with boundary $\Sigma := \partial M$.
A point $p \in \Sigma$ is \defn{singular} if $T_p\Sigma = H_p$.
We say that $\Sigma$ is \defn{nonsingular} if it contains no singular points.

Suppose now that $(M^3,T^{1,0},\theta)$ is a pseudohermitian three-manifold with nonsingular boundary $\Sigma := \partial M$.
Assume additionally that $\Sigma$ is oriented.
Then there is a unique $H \cap T\Sigma$-valued unit vector field $e_1$ such that $\nu := Je_1$ is inward-pointing.
Let $Z_1$ be the unique section of $T^{1,0}$ along $\Sigma$ such that $e_1 = \Real Z_1$.
A straightforward computation using~\eqref{eqn:ibp1} gives
\begin{equation}
 \label{eqn:ibp2}
 \int_M \langle \nablab u , \nablab w \rangle \, \theta \wedge d\theta = \int_M u\Delta_bw \, \theta \wedge d\theta - \oint_{\partial M} u\nabla_\nu w \, d\sigma ,
\end{equation}
where $d\sigma := \iota_{-\nu} (\theta \wedge d\theta)$.

Since the Tanaka--Webster connection preserves the contact form and the Levi form, $\nabla_{e_1} e_1$ is in the span of $e_2$.
The \defn{$p$-mean curvature} is the function $H$ defined by
\begin{equation*}
 \nabla_{e_1}e_1 := H\nu .
\end{equation*}

The \defn{CR Yamabe operator} $L^\theta \colon C^\infty(M) \to C^\infty(M)$ is
\begin{equation*}
 L^\theta u := \Delta_bu + \frac{R}{4}u .
\end{equation*}
The CR Yamabe operator is conformally covariant~\cite[Equation~(3.1)]{JerisonLee1987}:
If $\hat\theta = w^2\theta$, then
\begin{equation*}
 w^3 L^{\hat\theta} u = L^\theta(wu) .
\end{equation*}
When the contact form is clear from context, we shall write $L$ for $L^\theta$.
The \defn{CR Robin operator} $B^\theta \colon C^\infty(M) \to C^\infty(\Sigma)$ is
\begin{equation*}
 B^\theta u := -\nabla_\nu u + \frac{H}{3}u .
\end{equation*}
The conformal transformation law~\cite[Lemma~3.4]{JerisonLee1989} for the Tanaka--Webster connection implies that if $\hat\theta = w^2\theta$, then
\begin{equation*}
 w^2 B^{\hat\theta}u = B^\theta(wu) .
\end{equation*}
It follows that the \defn{CR Yamabe functional} $\mF \colon C^\infty(M) \to \R$,
\begin{equation*}
 \mF^\theta(u) := \int_M u \, Lu \, \theta \wedge d\theta + \oint_\Sigma u \, Bu \, d\sigma
\end{equation*}
is CR invariant;
indeed,
\begin{equation*}
 \mF^{w^2\theta}(u) = \mF^\theta(wu)
\end{equation*}
for all positive $w \in C^\infty(M)$ and all $u \in C^\infty(M)$.
Equation~\eqref{eqn:ibp2} implies that
\begin{equation*}
 \mF(u) = \int_M \left( \lvert \nablab u\rvert^2 + \frac{R}{4}u^2 \right) \, \theta \wedge d\theta + \frac{1}{3}\oint_\Sigma Hu^2 \, d\sigma .
\end{equation*}

\subsection{The standard sphere $\mathbb{S}^3$}
Let
\[
\theta = i \overline{\partial} (|\zeta|^2-1) = \frac{i}{2} (\zeta^1 d\zetabar^1 + \zeta^2 d\zetabar^2 - \zetabar^1 d\zeta^1 - \zetabar^2 d\zeta^2)
\]
be the contact form on $\mathbb{S}^3$. Denote 
\[
\theta^1 := \zeta^2 d\zeta^1 - \zeta^1 d\zeta^2.
\]
Then
\[
d\theta = i \theta^1 \wedge \theta^{\bar{1}}
\]
so $h_{1\bar{1}} = 1$.
Since
\[
d\theta^1 = -2i \theta^1 \wedge \theta,
\]
we have
\[
\omega_1{}^1 = -2i \theta.
\]
It follows that 
\[
d\omega_1{}^1 = 2 \theta^1 \wedge \theta^{\bar{1}},
\]
so the Tanaka--Webster scalar curvature is $R = 2$. The CR Yamabe operator is then
\begin{equation}
 \label{eqn:cr-yamabe-operator-sphere}
 L u := -(Z_1 Z_{\bar{1}} + Z_{\bar{1}} Z_1) u + \frac{1}{2} u
\end{equation}
where
\[
Z_1 := \overline{\zeta}^2 \frac{\partial}{\partial \zeta^1} - \overline{\zeta}^1 \frac{\partial}{\partial \zeta^2}.
\]
The fundamental solution to $L$ on $\mathbb{S}^3$ is~\cite[Theorem~2.1]{Geller1980}
\[
G(\zeta,\eta) := \frac{1}{8\pi} |1- \zeta \cdot \overline{\eta}|^{-1}.
\]

\subsection{The Clifford torus}

Recall the domain $\Omega \subset \mathbb{S}^3$ is given by
\[
\Omega = \{(\zeta^1,\zeta^2) \in \mathbb{S}^3 \colon |\zeta^1| < |\zeta^2|\}.
\]
Its boundary is the Clifford torus $\Sigma = \{(\zeta^1,\zeta^2) \colon |\zeta^1| = |\zeta^2| = \frac{1}{\sqrt{2}}\}$. 
Note that $\Sigma$ is contained in the open set $\{(\zeta^1, \zeta^2) \in \mathbb{S}^3 \colon |\zeta^1| \ne 0, |\zeta^2| \ne 0\}$. On this open set, define a frame of horizontal vectors
\begin{align*}
e_1 & := \frac{i}{\sqrt{2}} \left(\frac{\zeta^1 \zeta^2}{|\zeta^1 \zeta^2|} Z_1 - \frac{ \overline{\zeta}^1 \overline{\zeta}^2}{|\zeta^1 \zeta^2|} Z_{\bar{1}} \right) \\
& = - \sqrt{2} \, \Imag \left( |\zeta^2| \frac{\zeta^1}{|\zeta^1|} \frac{\partial}{\partial \zeta^1} - |\zeta^1| \frac{\zeta^2}{|\zeta^2|} \frac{\partial}{\partial \zeta^2} \right),
\end{align*}
and
\begin{align*}
\nu & := J e_1 = -\frac{1}{\sqrt{2}} \left(\frac{\zeta^1 \zeta^2}{|\zeta^1 \zeta^2|} Z_1 + \frac{ \overline{\zeta}^1 \overline{\zeta}^2}{|\zeta^1 \zeta^2|} Z_{\bar{1}} \right) \\
& = - \sqrt{2} \, \Real \left( |\zeta^2| \frac{\zeta^1}{|\zeta^1|} \frac{\partial}{\partial \zeta^1} - |\zeta^1| \frac{\zeta^2}{|\zeta^2|} \frac{\partial}{\partial \zeta^2} \right).
\end{align*}
The Reeb vector field on $\mathbb{S}^3$ is 
\begin{align*}
T &:= i \left( \zeta^1 \frac{\partial}{\partial \zeta^1} + \zeta^2 \frac{\partial}{\partial \zeta^2} - \overline{\zeta}^1 \frac{\partial}{\partial \overline{\zeta}^1} - \overline{\zeta}^2 \frac{\partial}{\partial \overline{\zeta}^2} \right) \\
&= -2 \Imag \left( \zeta^1 \frac{\partial}{\partial \zeta^1} + \zeta^2 \frac{\partial}{\partial \zeta^2} \right).
\end{align*}
For $t \in [0,\frac{\pi}{2\sqrt{2}})$, set 
\begin{equation*}
 \Omega_t = \left\{\zeta \in \mathbb{S}^3 \colon \frac{|\zeta^1|}{\cos(\frac{\pi}{4}+\frac{t}{\sqrt{2}})} < \frac{|\zeta^2|}{\sin(\frac{\pi}{4}+\frac{t}{\sqrt{2}}) } \right\}
\end{equation*} 
and 
$$\Sigma_t = \partial \Omega_t = \left\{\zeta \in \mathbb{S}^3 \colon |\zeta^1| = \cos\Bigl(\frac{\pi}{4}+\frac{t}{\sqrt{2}}\Bigr), |\zeta^2| = \sin\Bigl(\frac{\pi}{4}+\frac{t}{\sqrt{2}}\Bigr) \right\},$$ 
so that $\Omega_0 = \Omega$ and $\Sigma_0 = \Sigma$.
Then $\nu$ is the inward normal to $\partial \Omega_t$, and $e_1$ and $T$ are tangent to $\Sigma_t$.
In fact, a defining function for $\Omega_t$ is 
\[
\rho_t := \frac{|\zeta^1|^2}{\cos^2(\frac{\pi}{4}+\frac{t}{\sqrt{2}})} - \frac{|\zeta^2|^2}{\sin^2(\frac{\pi}{4}+\frac{t}{\sqrt{2}})},
\]
and the above expressions for $e_1$ and $T$ show that
\[
\nabla_{e_1} \rho_t = \nabla_T \rho_t = 0 \quad \text{on $\Sigma_t$} .
\]
Furthermore, 
\[
\nabla_{\nu} \rho_t = -\frac{\sqrt{2}}{\cos(\frac{\pi}{4}+\frac{t}{\sqrt{2}}) \sin(\frac{\pi}{4}+\frac{t}{\sqrt{2}})}  = -\frac{2\sqrt{2}}{\cos(\sqrt{2} t)} < 0
\]
so $\nu = Je_1$ is the inward-pointing normal to $\partial \Omega_t$. 
Moreover, we readily compute that $\nabla_{e_1}e_1=0$, and hence the Clifford torus is $p$-minimal.

\subsection{Fourier analysis} To parametrize $\Sigma_t$, it will be convenient to consider the lattice 
\[
\Lambda := 2\pi \{(m/2, n/2) \in \Z^2 \colon m \equiv n \pmod{2} \},
\]
and the associated abelian group $\R^2 / \Lambda$.
A fundamental domain for $\R^2 / \Lambda$ is given by $[-\pi,\pi) \times [-\pi/2,\pi/2)$; the identity element of $\R^2 / \Lambda$ will be written as $(0,0)$ using this identification. Any function defined on $[-\pi,\pi) \times [-\pi/2,\pi/2)$ can be lifted to $\R^2 / \Lambda$.
For instance, we define, for $(u,v) \in [-\pi,\pi) \times [-\pi/2,\pi/2)$, the natural norm for this problem, namely
\[
\norm{(u,v)} := \max\{|u|,|v|^{1/2}\},
\]
and it gives rise to a corresponding norm function on $\R^2 / \Lambda$.

The Lebesgue measure $du dv$ on $[-\pi,\pi) \times [-\pi/2,\pi/2)$ induces a Haar measure on $\R^2 / \Lambda$, with volume $2\pi^2$, and the dual lattice $(\frac{1}{2\pi}\Lambda)^*$ of $\frac{1}{2\pi}\Lambda$  is given by 
\[
(\frac{1}{2\pi}\Lambda)^* = \{(m,n) \in \Z^2 \colon m \equiv n \pmod{2} \}.
\]
An orthonormal basis of $L^2(\R^2 / \Lambda, dudv)$ is given by 
\[
\Big\{\frac{1}{\sqrt{2 \pi^2}} e^{i(mu+nv)} \colon (m,n) \in (\frac{1}{2\pi}\Lambda)^*  \Big\},
\] 
and hence Parseval's formula reads
\[
\iint_{\substack{|u| \leq \pi \\ |v| \leq \pi/2}} |F(u,v)|^2 du dv = \frac{1}{2\pi^2} \sum_{\substack{(m,n) \in \Z^2 \\ m \equiv n \pmod{2}}} |\widehat{F}(m,n)|^2,
\]
where 
\[
\widehat{F}(m,n) := \iint_{\substack{|u| \leq \pi \\ |v| \leq \pi/2}} F(u,v) e^{-i(mu+nv)} du dv.
\]
The convolution of two functions $F$ and $K$ on $\R^2 / \Lambda$ is given by any of the two equivalent expressions:
\[
\begin{split}
F*K(u',v') &:= \iint_{\substack{|u| \leq \pi \\ |v| \leq \pi/2}} F(u'-u,v'-v) K(u,v) du dv \\
&= \iint_{\substack{|u| \leq \pi \\ |v| \leq \pi/2}} F(u,v) K(u'-u,v'-v) du dv,
\end{split}
\]
for every $(u',v') \in \R^2 / \Lambda$. We have
\[
\widehat{F*K}(m,n) = \widehat{F}(m,n) \widehat{K}(m,n),
\]
thus 
\[
\iint_{\substack{|u| \leq \pi \\ |v| \leq \pi/2}} |F*K(u,v)|^2 du dv \leq \sup_{\substack{(m,n) \in \Z^2 \\ m \equiv n \pmod{2}}} |\widehat{K}(m,n)|^2 \iint_{\substack{|u| \leq \pi \\ |v| \leq \pi/2}} |F(u,v)|^2 du dv.
\]

\subsection{Parametrizing $\Sigma_t$}
\begin{figure}
\begin{tikzpicture}
\draw[->] (-2.5,0) -- (2.5,0);
\draw[->] (0,-2.5) -- (0,2.5);
\draw[dotted, fill=pink!40] (-1,1) -- (-1/2,3/2) -- (0,1) -- cycle;
\draw[dotted, fill=pink!40] (-1,-1) -- (-1/2,-1/2) -- (0,-1) -- cycle;
\draw[dotted, fill=green!20] (-1,1) -- (-3/2,1/2) -- (-1,0) -- cycle;
\draw[dotted, fill=green!20] (1,1) -- (1/2,1/2) -- (1,0) -- cycle;
\draw[dotted, fill=yellow!40] (1,-1) -- (3/2,-1/2) -- (1,0) -- cycle;
\draw[dotted, fill=yellow!40] (-1,-1) -- (-1/2,-1/2) -- (-1,0) -- cycle;
\draw[dotted, fill=blue!20] (0,-1) -- (1/2,-3/2) -- (1,-1) -- cycle;
\draw[dotted, fill=blue!20] (0,1) -- (1/2,1/2) -- (1,1) -- cycle;
\draw[dotted] (-1,-1) -- (-1,1) -- (1,1) -- (1,-1) -- cycle;
\draw[thick] (-3/2,1/2) -- (-1/2,3/2) -- (3/2,-1/2) -- (1/2,-3/2) -- cycle;
\end{tikzpicture}
\caption{The dotted square is $[-\pi,\pi]^2$. The tilted rectangle is the image of $(u,v) \in [-\pi,\pi] \times [-\pi/2,\pi/2] \mapsto (v+u,v-u) \in \R^2$. By translating the colored pieces in the tilted rectangle either vertically or horizontally by $2\pi$, we see that the tilted rectangle also parametrizes a copy of $(\R / 2\pi \Z)^2$. 
} 
\label{fig1}
\end{figure}
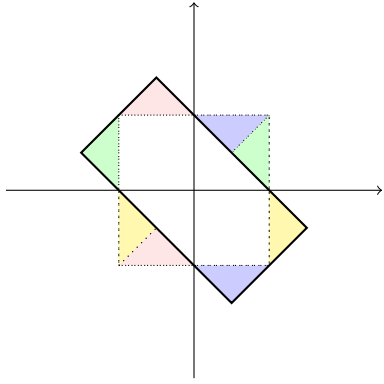

For $0 \leq t \leq \frac{\pi}{100}$, we will parametrize $\Sigma_t$ using 
\[
\begin{split}
(u,v) \in \R^2 / \Lambda 
\longmapsto  \Phi(t,u,v) := \left( \cos(\frac{\pi}{4}+\frac{t}{\sqrt{2}}) e^{i(v+u)},   \sin(\frac{\pi}{4}+\frac{t}{\sqrt{2}}) e^{i(v-u)} \right) \in \Sigma_t.
\end{split}
\]
One can see that the above map is a well-defined bijection by noting that $e^{i(v+u)} = e^{i(v-u)} = 1$ if and only if $(u,v) \in \Lambda$; alternatively, one observes that the map
\[
(u,v) \in [-\pi,\pi) \times [-\pi/2,\pi/2) \mapsto (v+u,v-u) \in (\R / 2\pi \Z)^2
\]
is a bijection (see Figure \ref{fig1}).
We have $\Phi_*(\frac{\partial}{\partial t}) = \nu$ and $\Phi_*(\frac{\partial}{\partial v}) = T$. 

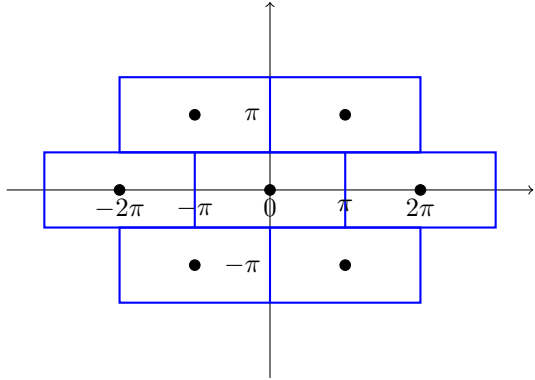
\begin{figure}
\begin{tikzpicture}
\draw[->] (-3.5,0) -- (3.5,0);
\draw[->] (0,-2.5) -- (0,2.5);
\draw[blue,thick] (-3,-0.5) -- (3,-0.5) -- (3,0.5) -- (-3,0.5) -- cycle;
\draw[blue,thick] (-1,-0.5) -- (-1,0.5);
\draw[blue,thick] (1,-0.5) -- (1,0.5);
\draw[blue,thick] (-2,0.5) -- (2,0.5) -- (2,1.5) -- (-2,1.5) -- cycle;
\draw[blue,thick] (0,0.5) -- (0,1.5);
\draw[blue,thick] (0,-0.5) -- (0,-1.5);
\draw[blue,thick] (-2,-0.5) -- (2,-0.5) -- (2,-1.5) -- (-2,-1.5) -- cycle;
\node[anchor=north] at (2,0) {$2\pi$};
\node[anchor=north] at (1,0) {$\pi$};
\node[anchor=north] at (0,0) {$0$};
\node[anchor=north] at (-1,0) {$-\pi$};
\node[anchor=north] at (-2,0) {$-2\pi$};
\node[anchor=east] at (0,1) {$\pi$};
\node[anchor=east] at (0,-1) {$-\pi$};
\filldraw[black] (0,0) circle (2pt);
\filldraw[black] (2,0) circle (2pt);
\filldraw[black] (-2,0) circle (2pt);
\filldraw[black] (1,1) circle (2pt);
\filldraw[black] (1,-1) circle (2pt);
\filldraw[black] (-1,1) circle (2pt);
\filldraw[black] (-1,-1) circle (2pt); 
\end{tikzpicture}
\caption{The dots represent the lattice $\Lambda$. Each blue rectangle is a fundamental domain for $\R^2 / \Lambda$. It can be identified with $\Sigma_t$ via the map $\Phi(t,\cdot,\cdot)$.
} 
\label{fig2}
\end{figure}

Under this parametrization of $\Sigma_t$, the volume element $d\sigma_t$ on $\Sigma_t$ is
\begin{align*}
d\sigma_t & =\iota_{-\nu} \theta \wedge d\theta \\
&= i \theta \wedge \iota_{\nu} (\theta^1 \wedge \theta^{\bar{1}}) \\
&= -\frac{1}{\sqrt{2}} i (\frac{\zeta^1 \zeta^2}{|\zeta^1 \zeta^2|} \theta \wedge \theta^{\bar{1}} - \frac{\overline{\zeta}^1 \overline{\zeta}^2}{|\zeta^1 \zeta^2|} \theta \wedge \theta^1) \\
&= -\frac{1}{\sqrt{2}} i (-\frac{\zeta^1 \zeta^2}{|\zeta^1 \zeta^2|} \frac{1}{2i} d\theta^{\bar{1}} - \frac{\overline{\zeta}^1 \overline{\zeta}^2}{|\zeta^1 \zeta^2|} \frac{1}{2i} d\theta^1) \\
&= -\frac{1}{\sqrt{2}} (\frac{\zeta^1 \zeta^2}{|\zeta^1 \zeta^2|} d\overline{\zeta}^1 \wedge d\overline{\zeta}^2 + \frac{\overline{\zeta}^1 \overline{\zeta}^2}{|\zeta^1 \zeta^2|} d\zeta^1 \wedge d\zeta^2) \\
&= -\sqrt{2} \cos(\frac{\pi}{4}+\frac{t}{\sqrt{2}}) \sin(\frac{\pi}{4}+\frac{t}{\sqrt{2}}) \Real \Big((-i (du + dv)) \wedge (-i (dv - du)) \Big) \\
&= \sqrt{2} \cos(\sqrt{2} t) du \wedge dv.
\end{align*}
In particular, setting $t=0$, the surface measure $d\sigma$ on $\Sigma$ is given by $\sqrt{2} du dv$.
If $f$ is an integrable function on $\Sigma_t$, we often write $F(u,v) := f(\Phi(0,u,v))$ so that 
\begin{equation} \label{eq:int_uv}
\oint_{\Sigma} f d\sigma = \sqrt{2} \iint_{\substack{|u| \leq \pi \\ |v| \leq \pi/2}} F(u,v) du dv.
\end{equation}

\section{The single layer potential} \label{sect3}
For $f \in C^1(\Sigma)$, the single layer potential is defined by
\[
\mathcal{S}f(\zeta) =  \int_{\eta \in \Sigma} f(\eta) G(\zeta,\eta) d\sigma(\eta), \quad \zeta \in \overline{\Omega}.
\]
Given $\zeta = (\zeta^1, \zeta^2) \in \Sigma$, the curve
\begin{equation} \label{eq:gamma_zeta}
\gamma_{\zeta}(t):=  \left( \sqrt{2} \cos(\frac{\pi}{4}+\frac{t}{\sqrt{2}}) \zeta^1 , \sqrt{2}  \sin(\frac{\pi}{4}+\frac{t}{\sqrt{2}}) \zeta^2 \right), \quad t \geq 0 ,
\end{equation}
satisfies $\gamma_{\zeta}(0) = \zeta$ and $\gamma_{\zeta}'(0) = \nu$ at $\zeta$. 

\begin{prop} \label{prop:derv_Sf}
For $f \in C^1(\Sigma)$, the one-sided normal derivative
\[
\nabla_{\nu} \mathcal{S}f(\zeta) := \lim_{t \to 0^+} \frac{\mathcal{S}f(\gamma_{\zeta}(t)) - \mathcal{S}f(\gamma_{\zeta}(0))}{t}
\]
exists at every $\zeta \in \Sigma$, and is given by
\[
\nabla_{\nu} \mathcal{S}f(\zeta) = (-\frac{1}{2} I + \mathcal{N})f(\zeta) ,
\]
where
\[
 \mathcal{N}f(\zeta) := p.v. \int_{\eta \in \Sigma} \langle \nu_{\zeta}, (\nabla_b)_{\zeta} G(\zeta,\eta) \rangle f(\eta) d\sigma(\eta), \quad \zeta \in \Sigma.
\]
\end{prop}

To prove this, first we understand the kernel $G(\zeta,\eta)$ of $\mathcal{S}$ in the $(u,v)$ coordinates. For $t \geq 0$ and $(u',v'), (u,v) \in \R^2 / \Lambda$, if $\zeta = \Phi(t,u',v')$ and $\eta = \Phi(0,u,v)$, then
\[
\begin{split}
& 1-\zeta \cdot \overline{\eta} \\
=& \,  1 - \cos (\frac{\pi}{4}) \cos(\frac{\pi}{4}+\frac{t}{\sqrt{2}}) e^{i(v'-v+u'-u)} - \sin(\frac{\pi}{4}) \sin(\frac{\pi}{4}+\frac{t}{\sqrt{2}}) e^{i(v'-v-u'+u)} \\
=& \, e^{i(v'-v)} \left[ e^{-i(v'-v)} - \cos(\frac{\pi}{4}) \cos(\frac{\pi}{4}+\frac{t}{\sqrt{2}}) e^{i(u'-u)} - \sin(\frac{\pi}{4}) \sin(\frac{\pi}{4}+\frac{t}{\sqrt{2}}) e^{-i(u'-u)} \right] \\
=& \, e^{i(v'-v)} \left[ \Big(\cos(v'-v)  - \cos(\frac{t}{\sqrt{2}}) \cos(u'-u) \Big) -i \Big(\sin(v'-v) - \sin(\frac{t}{\sqrt{2}}) \sin(u'-u) \Big) \right] 
\end{split}
\]
so from $G(\zeta,\eta) = \frac{1}{8\pi} \left| 1-\zeta \cdot \overline{\eta} \right|^{-1} $ we obtain
\begin{equation} \label{eq:G_to_k}
\begin{split}
G(\zeta,\eta) = \frac{1}{8\pi} k(\frac{t}{\sqrt{2}}, u'-u, v'-v)
\end{split}
\end{equation}
where for $t \geq 0$ and $(u,v) \in \R^2 / \Lambda$ we write
\[
k(t,u,v) := \Big[(\cos v - \cos t  \cos u)^2 + (\sin v - \sin t \sin u)^2 \Big]^{-1/2}.
\]

We collect a few facts we will need of $k(t,u,v)$:
\begin{Lem} \label{lem:Kt}
\begin{enumerate}[(a)]
    \item If $t \in (0,\pi/2)$, then $k(t,u,v)$ is a continuous function of $(u,v) \in \R^2 / \Lambda$, and so is $\frac{\partial}{\partial t} k(t,u,v)$. 
    \item If $t = 0$, then $k(t,u,v)$ is continuous at every $(u,v) \in \R^2 / \Lambda$ except at $(0,0)$, and is locally integrable near $(0,0)$.
    \item In fact, 
\[
|k(t,u,v)| \lesssim (|t| + \norm{(u,v)})^{-2}
\]
for all $t \in [0,\pi/4]$, $(u,v) \in \R^2 / \Lambda$, where the implicit constant is independent of $(t,u,v)$.
\item 
For all $(t,u,v) \in [0,\pi/4] \times (\R^2 / \Lambda) \setminus \{(0,0,0)\}$, we have
\[
\frac{\partial}{\partial t} k(t,u,v)
= - \frac{(\cos v - \cos t \cos u) \sin t \cos u - (\sin v - \sin t \sin u) \cos t \sin u}{ \Big[(\cos v - \cos t  \cos u)^2 + (\sin v - \sin t \sin u)^2 \Big]^{3/2}},
\]
and hence
\[
\Big| \frac{\partial}{\partial t} k(t,u,v) \Big| \lesssim (|t| + \norm{(u,v)})^{-3},
\]
where the implicit constant is independent of $(t,u,v)$.
\item  

We have
\begin{equation} \label{eq:int_Dtkt}
\lim_{t \to 0^+} \frac{1}{8\pi} \iint_{\substack{|u| \leq \pi \\ |v| \leq \pi/2}} \frac{\partial k}{\partial t}(t,u,v) du dv = -\frac{1}{2}.
\end{equation}

\end{enumerate}
\end{Lem}

We postpone the proof of this lemma to the end of this section.

\begin{proof}[Proof of Proposition~\ref{prop:derv_Sf}]
Using \eqref{eq:int_uv} and \eqref{eq:G_to_k}, if we write $F(u,v) := f(\Phi(0,u,v))$, we have
\begin{equation*} \label{eq:Sf_uv}
\mathcal{S}f(\zeta) = \frac{\sqrt{2}}{8\pi} \iint_{\substack{|u| \leq \pi \\ |v| \leq \pi/2}} F(u,v) k(\frac{t}{\sqrt{2}},u'-u,v'-v) du dv
\end{equation*}
whenever $\zeta = \Phi(t,u',v')$, $t \geq 0$. When $f \in L^1(\Sigma)$, the above integral converges absolutely at every point $\zeta$ by the estimate for $k(t,u,v)$ in Lemma~\ref{lem:Kt}. This integral is a convolution on the group $\R^2 / \Lambda$, and thus we can also write
\begin{equation} \label{eq:Sf_uv_nice}
\mathcal{S}f(\zeta) = \frac{\sqrt{2}}{8\pi} \iint_{\substack{|u| \leq \pi \\ |v| \leq \pi/2}} F(u'-u,v'-v) k(\frac{t}{\sqrt{2}},u,v) du dv.
\end{equation}
Similarly, when $f \in C^1(\Sigma)$ and $\zeta = \Phi(0,u',v')$, we have
\begin{align*} 
\mathcal{N}f(\zeta)
= \lim_{\varepsilon \to 0^+} \frac{1}{8\pi} \iint_{\substack{|u| \leq \pi, |v| \leq \pi/2 \\ \norm{(u'-u,v'-v)} \geq \varepsilon}} F(u,v) \frac{\partial k}{\partial t}(0,u'-u,v'-v) du dv \label{eq:DSf_uv} 
\end{align*}
is the convolution of $F$ with $\frac{1}{8\pi} p.v. \frac{\partial k}{\partial t}(0,\cdot)$ on $\R^2 / \Lambda$; the principal value exists when $f \in C^1(\Sigma)$ since $\frac{\partial k}{\partial t}(0,u,v)$ is odd in both $u$ and $v$, which implies 
\[
\iint_{\substack{|u| \leq \pi, |v| \leq \pi/2 \\ \norm{(u'-u,v'-v)} \geq \varepsilon}} \frac{\partial k}{\partial t}(0,u,v) du dv = 0
\] 
for every $\varepsilon > 0$, and $|\frac{\partial k}{\partial t}(0,u,v)| \lesssim \norm{(u,v)}^{-3}$.
In fact, one can write $\mathcal{N}f$ as an absolutely convergent integral:
\begin{equation} \label{eq:N_singular}
    \mathcal{N}f(\zeta)
    = \frac{1}{8\pi} \iint_{\substack{|u| \leq \pi \\ |v| \leq \pi/2}} \Big[F(u'-u,v'-v) - F(u',v')\Big] \frac{\partial k}{\partial t}(0,u,v) du dv.
\end{equation}

Now let $\zeta \in \Sigma$, and recall the curve $\gamma_{\zeta}(t)$ introduced in \eqref{eq:gamma_zeta}. Using \eqref{eq:Sf_uv_nice}, the bound for $k$ in Lemma~\ref{lem:Kt} and the Dominated Convergence Theorem, we see that
\[
\lim_{t \to 0^+} \mathcal{S}f(\gamma_{\zeta}(t)) = \mathcal{S}f(\gamma_{\zeta}(0)).
\]
Furthermore, from \eqref{eq:Sf_uv_nice} and our bounds for $\frac{\partial}{\partial t} k$, we see that for $t > 0$,
\[
\frac{d}{dt} \mathcal{S}f(\gamma_{\zeta}(t)) 
= \frac{1}{8\pi} \iint_{\substack{|u| \leq \pi \\ |v| \leq \pi/2}} F(u'-u,v'-v) \frac{\partial k}{\partial t}(\frac{t}{\sqrt{2}},u,v) du dv.
\]
Hence to prove Proposition~\ref{prop:derv_Sf}, by L'H{\^o}pital's rule and \eqref{eq:N_singular}, it suffices to show, under the assumption that $f \in C^1(\Sigma)$, that if $F(u,v) = f(\Phi(0,u,v))$ and $(u',v') \in [-\pi,\pi) \times [-\pi/2,\pi/2)$, then
\begin{equation} \label{eq:quant1}
\lim_{t \to 0^+} \frac{1}{8\pi} \iint_{\substack{|u| \leq \pi \\ |v| \leq \pi/2}} F(u'-u,v'-v) \frac{\partial k}{\partial t}(\frac{t}{\sqrt{2}},u,v) du dv
\end{equation}
exists and equals
\begin{equation} \label{eq:ptwise_goal}
-\frac{1}{2} F(u',v') + \frac{1}{8\pi} \iint_{\substack{|u| \leq \pi \\ |v| \leq \pi/2}} \Big[F(u'-u,v'-v) - F(u',v')\Big]\frac{\partial k}{\partial t}(0,u,v) du dv.
\end{equation}
To do so, we rewrite the expression inside the limit in \eqref{eq:quant1} as 
\begin{equation} \label{eq:ptwise_step1}
    \begin{split}
\frac{1}{8\pi} \iint_{\substack{|u| \leq \pi \\ |v| \leq \pi/2}} \Big[F(u'-u,v'-v) - F(u',v')\Big] \frac{\partial k}{\partial t}(t,u,v) du dv \\
+ \frac{1}{8\pi} F(u',v') \iint_{\substack{|u| \leq \pi \\ |v| \leq \pi/2}} \frac{\partial k}{\partial t}(t,u,v) du dv. 
\end{split}
\end{equation}
By \eqref{eq:int_Dtkt} in Lemma~\ref{lem:Kt}, as $t \to 0^+$, the second term in \eqref{eq:ptwise_step1} converges to $-\frac{1}{2} F(u',v')$, which appears as the first term of \eqref{eq:ptwise_goal}. Furthermore, since $f \in C^1(\Sigma)$, our earlier bound for  $\frac{\partial}{\partial t}k$ in Lemma~\ref{lem:Kt} and the Dominated Convergence Theorem shows that as $t \to 0^+$, the first term in \eqref{eq:ptwise_step1} converges to the second term in \eqref{eq:ptwise_goal}.
This completes the proof of Proposition~\ref{prop:derv_Sf}.
\end{proof}

\begin{proof}[Proof of Lemma~\ref{lem:Kt}]
First, for fixed $t$, the function $k(t,u,v)$ is doubly $2\pi$ periodic and invariant under $(u,v) \mapsto (u+\pi,v+\pi)$. Thus $k(t,u,v)$ is a well-defined function on $\R^2 / \Lambda$. 

To proceed further, we prove that 
\begin{multline*}
    (\cos v - \cos t \cos u)^2 + (\sin v - \sin t \sin u)^2 = 0 
    \Longleftrightarrow 
    \begin{cases}
    \cos t = 0 \\
    \cos u = 0 \\
    \sin v = \sin t \sin u
    \end{cases}
    \quad \text{or} \quad 
    \begin{cases}
    \sin t = 0 \\
    \sin u = 0 \\
    \cos v = \cos t \cos u
    \end{cases}.   
\end{multline*}
In fact, the former is equivalent to 
\[
    \begin{cases}
    \cos v = \cos t \cos u \\
    \sin v = \sin t \sin u,
    \end{cases}
\]
which implies, since $\cos^2 v +  \sin^2 v = 1$, that
\[
\begin{split}
1 
&= \cos^2 t \cos^2 u + \sin^2 t \sin^2 u \\
&= (\cos^2 t + \sin^2 t)(\cos^2 u + \sin^2 u) - \cos^2 t \sin^2 u - \sin^2 t \cos^2 u \\
&= 1 - \cos^2 t \sin^2 u - \sin^2 t \cos^2 u.
\end{split}
\]
As a result, $\cos t \sin u = \sin t \cos u = 0$. Either $\cos t = 0$, in which case $\sin t = \pm 1$ and $\cos u = 0$, or $\sin t = 0$, in which case $\cos t = \pm 1$ and $\sin u = 0$. From this the claimed equivalence follows.

If $t \in (0,\pi/2)$, then both $\cos t$ and $\sin t$ are non-zero, so the above implies
\[
    (\cos v - \cos t \cos u)^2 + (\sin v - \sin t \sin u)^2 > 0.
\]
This proves (a). If $t = 0$, then $\cos t = 1$, so the only zeroes of $(\cos v - \cos t \cos u)^2 + (\sin v - \sin t \sin u)^2$ are 
\[
\{(u,v) \in \R^2 \colon \sin u = 0, \cos v = \cos u\} = 2\pi \mathbb{Z}^2 \cup ((\pi,\pi)+2\pi\mathbb{Z}^2).
\]
This proves (b).

To prove (c), note that $k(t,u,v)$ is continuous on $[0,\pi/4] \times [-\pi,\pi] \times [-\pi/2,\pi/2] \setminus \{(0,0,0)\}$. So we only need to prove that
\[
(\cos v - \cos t \cos u)^2 + (\sin v - \sin t \sin u)^2 \gtrsim (|t|+|u|+|v|^{1/2})^4
\]
near $(t,u,v) = (0,0,0)$. But if we denote by
\begin{equation} \label{eq:g_def}
g(t,u,v) := (\cos v - \cos t  \cos u)^2 + (\sin v - \sin t \sin u)^2,
\end{equation} 
then by a Taylor expansion
\begin{equation} \label{eq:g_expansion}
\begin{split}
g(t,u,v) 
&= [(1-(1-\frac{t^2}{2})(1-\frac{u^2}{2}) + O^{5/2}]^2 + (v - tu + O^{5/2})^2 \\
&= (\frac{t^2 + u^2}{2})^2 + (v-tu)^2 + O^{9/2}\\
&= ((\frac{t^2 + u^2}{2})^2 + (v-tu)^2 ) (1 + O^{1/2})
\end{split}
\end{equation}
where $O^{j}$ is an error bounded by $C_j (|t| + |u| + |v|^{1/2})^{j}$ for some absolute constant $C_j$. Our claim now follows from the fact that
\begin{equation} \label{eq:norm-equiv}
(\frac{t^2 + u^2}{2})^2 + (v-tu)^2 \simeq (|t| + |u| + |v|^{1/2})^4
\end{equation}
uniformly in $t, u, v$ (which can be verified by noting that the left hand side is a positive continuous function on the set where $|t|+|u|+|v|^{1/2} = 1$, and appealing to homogeneity).


(d) follows from
\[
\begin{split}
|\frac{\partial}{\partial t} k(t,u,v)| &\leq \frac{|\sin t \cos u|+ |\cos t \sin u|}{(\cos v - \cos t  \cos u)^2 + (\sin v - \sin t \sin u)^2} \leq \frac{|t|+|u|}{g(t,u,v)}
\end{split}
\]
and the lower bound $g(t,u,v) \gtrsim |t|+|u|+|v|^{1/2}$ on $[0,\pi/4] \times [-\pi,\pi] \times [-\pi/2,\pi/2]$ that we proved above.

Finally, we prove \eqref{eq:int_Dtkt} in (e). 
Since $\frac{\partial k}{\partial t}(t,u,v)$ converges uniformly to $0$ on $[-\pi,\pi] \times [-\pi/2,\pi/2] \setminus [-\varepsilon,\varepsilon] \times [-\varepsilon^2,\varepsilon^2]$ for any $\varepsilon > 0$, it suffices to show that there exists $\varepsilon > 0$ such that
\begin{equation*}
    \lim_{t \to 0^+} \iint_{\substack{|u| \leq \varepsilon \\ |v| \leq \varepsilon^2}} \frac{\partial k}{\partial t}(t,u,v) du dv = -4\pi.
\end{equation*}
From 
\[
\frac{\partial}{\partial t} k(t,u,v)
= - \frac{(\cos v - \cos t \cos u) \sin t \cos u - (\sin v - \sin t \sin u) \cos t \sin u}{ g(t,u,v)^{3/2}},
\]
where $g(t,u,v)$ was defined in \eqref{eq:g_def}, it suffices to show that
\begin{equation*}
    \lim_{t \to 0^+} \iint_{\substack{|u| \leq \varepsilon \\ |v| \leq \varepsilon^2}} \frac{\sin t \cos u (\cos v - \cos t \cos u) }{g(t,u,v)^{3/2}} du dv = 4\pi,
\end{equation*}
\begin{equation*}
\lim_{t \to 0^+} \iint_{\substack{|u| \leq \varepsilon \\ |v| \leq \varepsilon^2}} \frac{\sin u \sin v}{g(t,u,v)^{3/2}} du dv = 4\pi
\end{equation*}
and
\begin{equation*}
\lim_{t \to 0^+} \iint_{\substack{|u| \leq \varepsilon \\ |v| \leq \varepsilon^2}} \frac{\sin t \sin^2 u}{g(t,u,v)^{3/2}} du dv = 4\pi. 
\end{equation*}
We do so by Taylor expanding the numerators and denominators of the integrands.
By choosing $\varepsilon > 0$ sufficiently small, we have, whenever $t \in [0,\varepsilon]$, $|u| \leq \varepsilon$ and $|v| \leq \varepsilon^2$, that $(1+O^{1/2})^{-3/2} = 1 + O^{1/2}$, and hence from \eqref{eq:g_expansion}
\begin{align*}
g(t,u,v)^{-3/2} = [ (\frac{t^2 + u^2}{2})^2 + (v-tu)^2 ]^{-3/2} (1 + O^{1/2}).
\end{align*}
Similarly, for $t \in [0,\varepsilon]$, $|u| \leq \varepsilon$ and $|v| \leq \varepsilon^2$, we have
\[
\sin t \cos u (\cos v - \cos t \cos u) = t \Big[ \frac{( t^2 + u^2)}{2} + O^{5/2} \Big]
\]
and 
\[
\sin t \sin^2 u = t \Big [ u^2 + O^{5/2} \Big].
\]
Thus 
\begin{align*}
    &\lim_{t \to 0^+} \iint_{\substack{|u| \leq \varepsilon \\ |v| \leq \varepsilon^2}} \frac{\sin t \cos u (\cos v - \cos t \cos u) }{g(t,u,v)^{3/2}} du dv  \\
    &= \lim_{t \to 0^+} \iint_{\substack{|u| \leq \varepsilon \\ |v| \leq \varepsilon^2}} \frac{t \Big[\frac{ (t^2 + u^2)}{2} + O^{5/2} \Big]}{[ (\frac{t^2 + u^2}{2})^2 + (v-tu)^2 ]^{3/2}} du dv. 
\end{align*}
The error term involving $O^{5/2}$ can be bounded using 
\eqref{eq:norm-equiv}.
Thus
\begin{equation} \label{eq:Taylor_error}
\begin{split}
& \iint_{\substack{|u| \leq \varepsilon \\ |v| \leq \varepsilon^2}} \frac{t O^{5/2} }{[ (\frac{t^2 + u^2}{2})^2 + (v-tu)^2 ]^{3/2}} du dv \\
& \leq C t \iint_{\substack{|u| \leq \varepsilon \\ |v| \leq \varepsilon^2}} (|t|+|u|+|v|^{1/2})^{-7/2} du dv \\
&= C t \iint_{\substack{|u| \leq \varepsilon/t \\ |v| \leq \varepsilon^2/t^2}} (t+|tu|+|t^2 v|^{1/2})^{-7/2} t^3 du dv \\
&= C t^{1/2} \iint_{\substack{|u| \leq \varepsilon/t \\ |v| \leq \varepsilon^2/t^2}} (1+|u|+| v|^{1/2})^{-7/2} du dv
\end{split}
\end{equation}
where we performed a change of variables $(u,v) \mapsto (ut,vt^2)$ in the second-to-last line, and the last line is $O(t^{1/2}$) as $t \to 0^+$ since the double integral converges to a finite limit (we have $\iint_{\mathbb{R}^2} (1+|u|+|v|^{1/2})^{\alpha} du dv < \infty$ if and only if $\alpha < -3$). Furthermore, the main term in our earlier computation can be identified, via the same change of variable, as
\begin{align*}
    & \lim_{t \to 0^+} \iint_{\substack{|u| \leq \varepsilon/t \\ |v| \leq \varepsilon^2/t^2}} \frac{t \frac{t^2+(ut)^2}{2}}{[ (\frac{t^2 + (ut)^2}{2})^2 + (v t^2 - t(ut))^2 ]^{3/2}} t^3 du dv \\
    &= \iint_{\mathbb{R}^2} \frac{\frac{1+u^2}{2}}{[ (\frac{1+u^2}{2})^2 + (v-u)^2 ]^{3/2}} du dv = 4\pi.
\end{align*}
Similarly,
\begin{align*}
    &\lim_{t \to 0^+} \iint_{\substack{|u| \leq \varepsilon \\ |v| \leq \varepsilon^2}} \frac{\sin t \sin^2 u}{g(t,u,v)^{3/2}} du dv \\
    &= \lim_{t \to 0^+} \iint_{\substack{|u| \leq \varepsilon \\ |v| \leq \varepsilon^2}} \frac{t \Big[u^2 + O^{5/2} \Big]}{[ (\frac{t^2 + u^2}{2})^2 + (v-tu)^2 ]^{3/2}} du dv \\
    &= \iint_{\mathbb{R}^2} \frac{u^2}{[ (\frac{1 + u^2}{2})^2 + (v-u)^2 ]^{3/2}} du dv = 4\pi.
\end{align*}
Finally, use $\int_{|v| \leq \varepsilon^2} f(v) dv = \frac{1}{2} \int_{|v| \leq \varepsilon^2} [f(v) + f(-v)] dv$ and write
\begin{equation} \label{eq:g_diff}
\begin{split}
    &\iint_{\substack{|u| \leq \varepsilon \\ |v| \leq \varepsilon^2}} \frac{\sin u \sin v}{g(t,u,v)^{3/2}} du dv \\
    &= \frac{1}{2} \iint_{\substack{|u| \leq \varepsilon \\ |v| \leq \varepsilon^2}} \frac{\sin u \sin v}{g(t,u,v)^{3/2}} - \frac{\sin u \sin v}{g(t,u,-v)^{3/2}}du dv \\
    &= \frac{1}{2} \iint_{\substack{|u| \leq \varepsilon \\ |v| \leq \varepsilon^2}} \frac{\sin u \sin v [g(t,u,-v)-g(t,u,v)][g(t,u,-v)^2 + g(t,u,-v) g(t,u,v) + g(t,u,v)^2]}{g(t,u,v)^{3/2} g(t,u,-v)^{3/2} [g(t,u,-v)^{3/2}+g(t,u,v)^{3/2}]} du dv.
\end{split}
\end{equation}
We Taylor expand the numerator by noting
\[
\sin u \sin v = u v + O^{7/2},
\]
and
\[
\begin{split}
&g(t,u,-v)^2 + g(t,u,-v) g(t,u,v) + g(t,u,v)^2 \\
&= h(t,u,-v)^2 + h(t,u,-v) h(t,u,v) + h(t,u,v)^2 + O^{17/2}
\end{split}
\]
where we abbreviated 
\[
h(t,u,v) := (\frac{t^2 + u^2}{2})^2 + (v-tu)^2,
\]
and 
\[
\begin{split}
g(t,u,-v) - g(t,u,v) &= 4 \sin t \sin u \sin v \\
&= 4 t u v + t O^{7/2} \\
&= h(t,u,-v) - h(t,u,v) + t O^{7/2}.
\end{split}
\]
We also Taylor expand the denominator, yielding
\begin{align*}
&g(t,u,v)^{3/2} g(t,u,-v)^{3/2} [g(t,u,-v)^{3/2}+g(t,u,v)^{3/2}] \\
&= h(t,u,v)^{3/2} h(t,u,-v)^{3/2} [h(t,u,-v)^{3/2}+h(t,u,v)^{3/2}] (1+O^{1/2}).
\end{align*}
Thus
\begin{align*}
    &\lim_{t\to 0^+} \iint_{\substack{|u| \leq \varepsilon \\ |v| \leq \varepsilon^2}} \frac{\sin u \sin v}{g(t,u,v)^{3/2}} du dv \\
    &= \lim_{t\to 0^+} \frac{1}{2} \iint_{\substack{|u| \leq \varepsilon \\ |v| \leq \varepsilon^2}} \frac{u v [h(t,u,-v)-h(t,u,v)][h(t,u,-v)^2 + h(t,u,-v) h(t,u,v) + h(t,u,v)^2]}{h(t,u,v)^{3/2} h(t,u,-v)^{3/2} [h(t,u,-v)^{3/2}+h(t,u,v)^{3/2}]} du dv \\
    &\qquad + \lim_{t\to 0^+} \iint_{\substack{|u| \leq \varepsilon \\ |v| \leq \varepsilon^2}} \frac{t O^{29/2}}{h(t,u,v)^{3/2} h(t,u,-v)^{3/2} [h(t,u,-v)^{3/2}+h(t,u,v)^{3/2}]} du dv.
\end{align*}
The second limit is zero by the same argument of the proof of \eqref{eq:Taylor_error}. In the first limit, we are considering the limit as $t \to 0^+$ of a double integral reminiscient of the last expression of \eqref{eq:g_diff}. By reversing the derivation of \eqref{eq:g_diff}, except that we write $h$ in place of $g$, we see that the first limit is equal to
\begin{align*}
    & \lim_{t\to 0^+} \iint_{\substack{|u| \leq \varepsilon \\ |v| \leq \varepsilon^2}} \frac{u v}{h(t,u,v)^{3/2}} du dv \\
    &= \lim_{t\to 0^+} \iint_{\substack{|u| \leq \varepsilon \\ |v+tu| \leq \varepsilon^2}} \frac{(v+tu) u}{[(\frac{t^2 + u^2}{2})^2 + v^2]^{3/2}} du dv \\
    &= \lim_{t\to 0^+} \iint_{\substack{ |u| \leq \varepsilon \\ -tu-\varepsilon^2 \leq v \leq -tu+\varepsilon^2}} \frac{u v}{[(\frac{t^2 + u^2}{2})^2 + v^2]^{3/2}} du dv + \lim_{t\to 0^+} \iint_{\substack{|u| \leq \varepsilon \\ |v+tu| \leq \varepsilon^2}}  \frac{t u^2}{[(\frac{t^2 + u^2}{2})^2 + v^2]^{3/2}} du dv \\
    &= I + II
\end{align*}
We have 
\begin{align*}
I=&\lim_{t \to 0^+} \int_{|u| \leq \varepsilon} \frac{u}{[(\frac{t^2+u^2}{2})^2 + (-tu-\varepsilon^2)^2]^{1/2}} - \frac{u}{[(\frac{t^2+u^2}{2})^2 + (-tu+\varepsilon^2)^2]^{1/2}} du \\
=& \int_{|u| \leq \varepsilon} \frac{u}{(\frac{u^2}{4}+\varepsilon^4)^{1/2}} du - \int_{|u| \leq \varepsilon} \frac{u}{(\frac{u^2}{4}+\varepsilon^4)^{1/2}} du = 0 - 0 = 0
\end{align*}
and using a change of variable $(u,v) \mapsto (tu, t^2 v)$,
\[
II = \lim_{t\to 0^+} \iint_{\substack{|u| \leq \frac{\varepsilon}{t} \\ |v+u| \leq \frac{\varepsilon^2}{t^2}}}  \frac{u^2}{[(\frac{1 + u^2}{2})^2 + v^2]^{3/2}} du dv  = \iint_{\mathbb{R}^2} \frac{u^2}{[(\frac{1+u^2}{2})^2+v^2]^{3/2}} dudv = 4\pi.
\]
Thus 
\[
\lim_{t\to 0^+} \iint_{\substack{|u| \leq \varepsilon \\ |v| \leq \varepsilon^2}} \frac{\sin u \sin v}{g(t,u,v)^{3/2}} du dv = 4\pi,
\]
as desired.
\end{proof}

\section{The classical solution to the Neumann boundary value problem} \label{sect4}

In this section we prove Theorem~\ref{thm:Neumann}.

\subsection{Existence} We are now going to invert $-\frac{1}{2}I + \mathcal{N}$ on $L^2$, which will give us the solution to the Neumann boundary value problem for the CR Yamabe operator when the boundary data $h$ is in $C^{\infty}(\Sigma)$.

\begin{prop} \label{prop:4.1}
We have $\|\mathcal{N}\|_{L^2 \to L^2} < \frac{1}{2}$. As a result, for any $h \in L^2(\Sigma)$, there exists a unique $f \in L^2(\Sigma)$ such that 
\[
(-\frac{1}{2} I + \mathcal{N})f = h.
\]
Furthermore, if $h \in C^{\infty}(\Sigma)$, then $f \in C^{\infty}(\Sigma)$. 
\end{prop}

The main thrust is in showing that $\|\mathcal{N}\|_{L^2 \to L^2} < \frac{1}{2}$. Once that is proved, then $-\frac{1}{2}I + \mathcal{N}$ is invertible on $L^2$, and is given by a multiplier operator whose coefficients are bounded below. Hence $(-\frac{1}{2} I + \mathcal{N})^{-1}$ maps $W^{k,2}(\Sigma)$ to $W^{k,2}(\Sigma)$ for any $k \geq 0$, and hence maps $C^{\infty}(\Sigma)$ to $C^{\infty}(\Sigma)$. As a result, we obtain the existence assertion in Theorem~\ref{thm:Neumann}.

Below we prove that $\|\mathcal{N}\|_{L^2 \to L^2} < \frac{1}{2}$.

Indeed, since $\mathcal{N}f(\Phi(0,u',v'))$ is the convolution of $f(\Phi(0,u,v))$ with $\frac{1}{8\pi} K(u,v)$ where $K(u,v) := p.v. \frac{\partial k}{\partial t}(0,u,v)$ on $\R^2 / \Lambda$, we have
\[
\begin{split}
\int_{\Sigma} |\mathcal{N} f|^2 d\sigma 
&= \sqrt{2} \iint_{\substack{|u| \leq \pi \\ |v| \leq \pi/2}} |\mathcal{N}f(\Phi(0,u,v))|^2 du dv \\
&\leq \sup_{\substack{(m,n) \in \Z^2 \\ m \equiv n \pmod{2}}} |\frac{1}{8\pi} \widehat{K}(m,n)|^2 \sqrt{2} \iint_{\substack{|u| \leq \pi \\ |v| \leq \pi/2}} |f(\Phi(0,u,v))|^2 du dv \\
&= \sup_{\substack{(m,n) \in \Z^2 \\ m \equiv n \pmod{2}}} |\frac{1}{8\pi} \widehat{K}(m,n)|^2 \int_{\Sigma} |f|^2 d\sigma ,
\end{split}
\]
where
\[
\begin{split}
\widehat{K}(m,n) &=  \lim_{\varepsilon \to 0^+} \iint_{\substack{|u| \leq \pi, |v| \leq \pi/2 \\ \norm{(u,v)} \geq \varepsilon}} \frac{\partial k}{\partial t}(0,u,v) e^{-i(mu+nv)} du dv.
\end{split}
\]
Thus we need to show that
\begin{equation} \label{eq:K_multiplier_bdd_new}
\sup_{\substack{(m,n) \in \Z^2 \\ m \equiv n \pmod{2}}} |\widehat{K}(m,n)| < \frac{8\pi}{2} = 4\pi \simeq 12.56.
\end{equation}
Since $\frac{\partial k}{\partial t}(0,u,v)$ is odd in both $u$ and $v$, we have the following expression for $\widehat{K}(m,n)$:
\begin{equation} \label{eq:K_hat_convenient}
\begin{split}
\widehat{K}(m,n) &= -\iint_{\substack{|u| \leq \pi \\ |v| \leq \pi/2}} K(u,v) \sin(mu) \sin(nv) du dv.
\end{split}
\end{equation}

To bound this integral and establish \eqref{eq:K_multiplier_bdd_new}, it will be convenient to approximate $K(u,v)$ by a Taylor expanding its numerator and denominator. 
Let's write
\[
K(u,v) = \frac{n(u,v)}{d(u,v)^6}
\]
where
\[
n(u,v) := \sin u \sin v, \qquad d(u,v) := [(\cos u - \cos v)^2 + \sin^2 v]^{1/4}.
\]
We approximate $n(u,v)$ and $d(u,v)$ by
\[
n_0(u,v) := uv \qquad \text{and} \qquad d_0(u,v) := \Big[\frac{u^4}{4}+v^2 \Big]^{1/4}
\]
respectively, so that the homogeneous function
\[
K_{\circ}(u,v) := \frac{n_0(u,v)}{d_0(u,v)^6}
\]
will be a good approximation of $K(u,v)$ near $(u,v) = (0,0)$, along with its derivatives; this homogeneous function can be integrated easily near the origin, so in \eqref{eq:K_hat_convenient} we can approximate $K(u,v)$ by $K_{\circ}(u,v)$, and  eventually obtain \eqref{eq:K_multiplier_bdd_new}. In fact, it will be convenient to note that
\begin{equation} \label{eq:d_lower}
d(u,v)^4
\geq  d_0(u,v)^4 - \Big(\frac{u^6}{4!} + \frac{u^2v^2}{2} + \frac{2 u^8}{8!} + \frac{2 v^4}{4!} \Big)
\end{equation}
for all $(u,v) \in \R^2$, which follows since
\[
d(u,v)^4 = \cos^2 u - 2 \cos u \cos v + 1 = \frac{1}{2}(1+\cos 2u) - 2 \cos u + 1 + 2 \cos u (1 - \cos v)
\]
whereas
\begin{align*}
\frac{1}{2}(1+\cos 2u) - 2 \cos u + 1 
&\geq \frac{1}{2} \left(1+1-\frac{(2u)^2}{2!} + \frac{(2u)^4}{4!} - \frac{(2u)^6}{6!} \right) - 2 \left( 1 - \frac{u^2}{2!} + \frac{u^4}{4!} - \frac{u^6}{6!} + \frac{u^8}{8!} \right) + 1 \\
&= \frac{u^4}{4}-\frac{u^6}{4!} - \frac{2 u^8}{8!} 
\end{align*}
and
\begin{align*}
2 \cos u (1 - \cos v) &= \int_0^v \int_0^w 2 \cos u \cos x \, dx dw \\
&= \int_0^v \int_0^w \cos(u+x)+\cos(u-x) \, dx dw \\
&\geq \int_0^v \int_0^w \left(1-\frac{(u+x)^2}{2} \right) + \left(1-\frac{(u-x)^2}{2} \right) \, dx dw \\
&= \int_0^v \int_0^w (2-u^2-x^2) \, dx dw \\
&= v^2 - \frac{u^2 v^2}{2} - \frac{2 v^4}{4!}.
\end{align*}

To carry out the strategy mentioned above, we differentiate $K(u,v)$ and write
\[
\frac{\partial K}{\partial u} = \frac{n_u(u,v)}{d(u,v)^{10}}, \qquad \frac{\partial K}{\partial v} = \frac{n_v(u,v)}{d(u,v)^{10}}, \qquad
\frac{\partial^2 K}{\partial u^2} = \frac{n_{uu}(u,v)}{d(u,v)^{14}},
\]
where we have, via the half-angle formulae, the following expressions for the numerators:
\begin{align*}
 n_u(u,v) & := \sin(v) \left( -20 \sin^4(\frac{u}{2}) + 4 \sin^2(\frac{v}{2}) \right) \\
 & \quad + \sin(v) \left( 16 \sin^6(\frac{u}{2}) + 2 \sin^2(u) \sin^2(\frac{v}{2}) \right) , \\
 n_v(u,v) & := \sin(u) \left(4 \sin^4(\frac{u}{2}) -8 \sin^2(\frac{v}{2}) \right) \\
 & \quad + \sin(u) \sin^2(\frac{v}{2}) \left[ 8 \sin^2(\frac{u}{2}) \Big(1+\cos^2(\frac{u}{2}) \Big)+4 \cos(u) \sin^2(\frac{v}{2}) \right] 
\end{align*}
\begin{align*}
n_{uu}(u,v) & :=  120 \sin(u) \sin^2(\frac{u}{2}) \sin(v) (\sin^4(\frac{u}{2}) - \sin^2(\frac{v}{2}))  \\
 & \quad - 64 \sin(u) \sin^8(\frac{u}{2}) \sin(v) \\
 & \quad - 8 \Big(\cos(\frac{u}{2}) + \cos(\frac{3u}{2})\Big) \sin^5(\frac{u}{2}) \sin(v) \sin^2(\frac{v}{2}) \\
 & \quad + 4 \Big(14 + \sin^2(u)\Big) \sin(u) \sin(v) \sin^4(\frac{v}{2}).
\end{align*}
Similarly we write
\[
 \frac{\partial K_{\circ}}{\partial u} = \frac{n_{0u}(u,v)}{d_0(u,v)^{10}}, \qquad \frac{\partial K_{\circ}}{\partial v}  = \frac{n_{0v}(u,v)}{d_0(u,v)^{10}}, \qquad
 \frac{\partial^2 K_{\circ}}{\partial u^2} = \frac{n_{0uu}(u,v)}{d_0(u,v)^{14}},
\]
where
\[
n_{0u}(u,v) := v (-\frac{5}{4} u^4 + v^2), \qquad
n_{0v}(u,v) := u (\frac{u^4}{4} - 2 v^2) , \qquad
n_{0uu}(u,v) := \frac{15}{2} u^3 v (\frac{u^4}{4} - v^2).
\]
Later we will use 
\begin{equation} \label{eq:nu_est1}
|n_u(u,v)| \leq |n_{0u}(u,v)| + \frac{1451 |u^6 v| }{3072} + \frac{|u^2 v^3|}{2} +  \frac{v^4}{3} 
\end{equation}
\begin{equation} \label{eq:nv_est1}
|n_v(u,v)| \leq |n_{0v}(u,v)| + \frac{683|u^7| }{15360}  + |u^3 v^2| + \frac{7|u v^3|}{6} 
\end{equation}
\begin{equation} \label{eq:nuu_est1}
|n_{uu}(u,v)| \leq |n_{0uu}(u,v)| + \frac{1195 |u^9 v|}{2048} + \frac{|u^5 v^3|}{8} + \frac{5|u^3 v^4|}{2}  + \frac{15 |u v^5|}{4}
\end{equation}
In fact, we have
\[
\left| \frac{d^6}{du^6} \sin^4(\frac{u}{2}) \right| = \left| \frac{\cos(u)}{2}  - 8 \cos (2u) \right| \leq \frac{2049}{256} \quad \text{for all $u \in \R$}
\]
and
\[
\left| \frac{d^3}{dv^3} \sin^2(\frac{v}{2}) \right| = \left|-\frac{\sin v}{2} \right| \leq \frac{1}{2} \quad \text{for all $v \in \R$}.
\]
Thus Taylor expansion (up to order 5 and 3 respectively) gives
\[
\left| \sin^4(\frac{u}{2})  - \frac{u^4}{16} \right| \leq \frac{2049}{256} \frac{u^6}{6!}
\]
and
\[
\left| \sin^2(\frac{v}{2}) - \frac{v^2}{4} \right| \leq \frac{1}{2} \frac{|v|^3}{3!}.
\]
As a result, we can estimate the main terms for the expressions for $n_u$, $n_v$ and $n_{uu}$:
\[
\begin{split}
\left|\sin(v) \left( -20 \sin^4(\frac{u}{2}) + 4 \sin^2(\frac{v}{2}) \right) \right| & \leq |v| \left( \left| -20 \cdot \frac{u^4}{16} + 4 \cdot \frac{v^2}{4} \right| + 20 \cdot \frac{2049}{256} \frac{u^6}{6!} + 4 \cdot \frac{1}{2} \frac{|v|^3}{3!} \right) \\
&= |n_{0u}(u,v)| +  \frac{683}{3072} |u^6 v| + \frac{1}{3} v^4, 
\end{split}
\]
\[
\begin{split}
\left| \sin(u) \left( 4 \sin^4(\frac{u}{2}) - 8 \sin^2(\frac{v}{2}) \right) \right| & \leq |u| \left( \left| 4 \cdot \frac{u^4}{16} - 8 \cdot \frac{v^2}{4} \right| + 4 \cdot \frac{2049}{256} \frac{u^6}{6!} + 8 \cdot \frac{1}{2} \frac{|v|^3}{3!} \right) \\
&= |n_{0v}(u,v)| + \frac{683}{15360} |u^7| + \frac{2}{3} |u v^3|,
\end{split}
\]
and
\[
\begin{split}
\left|120 \sin(u) \sin^2(\frac{u}{2}) \sin(v) \left( \sin^4(\frac{u}{2}) - \sin^2(\frac{v}{2}) \right) \right| & \leq 30 |u|^3 |v| \left( \left| \frac{u^4}{16} - \frac{v^2}{4} \right| + \frac{2049}{256} \frac{u^6}{6!} + \frac{1}{2} \frac{|v|^3}{3!} \right) \\
&= |n_{0uu}(u,v)| + \frac{683}{2048} |u^9 v| + \frac{5}{2} |u^3 v^4|.
\end{split}
\]
Furthermore, the error terms can be bounded using
\[
\begin{split}
\left| \sin(v) \left( 16 \sin^6(\frac{u}{2}) + 2 \sin^2(u) \sin^2(\frac{v}{2}) \right) \right|
& \leq |v| \left( 16 \frac{|u|^6}{2^6} + 2 |u|^2 \frac{|v|^2}{2^2} \right) \\
& = \frac{1}{4} |u^6 v| + \frac{1}{2} |u^2 v^3|,
\end{split}
\]
\[
\begin{split}
\left| \sin(u) \sin^2(\frac{v}{2}) \left[ 8 \sin^2(\frac{u}{2}) \Big(1+\cos^2(\frac{u}{2}) \Big)+4 \cos(u) \sin^2(\frac{v}{2}) \right] \right|
& \leq |\frac{u v^2}{4}| \left[ 8 \frac{u^2}{2^2} (1+1) + 4 \cdot 1 \cdot 1 \cdot \frac{|v|}{2} \right] \\
& = |u^3 v^2| + \frac{1}{2} |u v^3|,
\end{split}
\]
\[
\left| - 64 \sin(u) \sin^8(\frac{u}{2}) \sin(v)  \right|
\leq 64 |u| \frac{u^8}{2^8} |v| = \frac{1}{4}{|u^9 v|},
\]
\[
\left|  - 8 \Big(\cos(\frac{u}{2}) + \cos(\frac{3u}{2})\Big) \sin^5(\frac{u}{2}) \sin(v) \sin^2(\frac{v}{2}) \\
  \right|
\leq 8 (1 + 1) \frac{|u|^5}{2^5} |v| \frac{|v|^2}{2^2} = \frac{1}{8}|u^5 v^3|,
\]
and
\[
\left| 4 \Big(14 + \sin^2(u)\Big) \sin(u) \sin(v) \sin^4(\frac{v}{2}) \right| \leq 4(14+1)|u v| \frac{v^4}{2^4} = \frac{15}{4} |u v^5|.
\]
The above nine inequalities together imply \eqref{eq:nu_est1}, \eqref{eq:nv_est1} and \eqref{eq:nuu_est1}.

Now let 
\[
b := \frac{2}{5}.
\]
For $0 < a \leq 5$, let 
\[
e_0(a) := \sup\left\{\frac{u^6}{4!} + \frac{u^2v^2}{2} +  \frac{2 u^6 (ab)^2}{8!} + \frac{2 |v|^3 b^2}{4!} \colon d_0(u,v) = 1 \right\}
\]
and
\[
e_1(a):= 1 - e_0(a) d_0(ab,b^2)^{2}.
\]
By homogeneity, we then have
\begin{equation} \label{eq:e0use}
\frac{u^6}{4!} + \frac{u^2v^2}{2} +  \frac{2 u^6 (ab)^2}{8!} + \frac{2 |v|^3 b^2}{4!} \leq e_0(a) d_0(u,v)^{6} \quad \text{for all $(u,v) \in \R^2$}
\end{equation}
Furthermore, $e_1(a)$ is a decreasing function of $a \in [0,\infty)$, with $e_1(5) > 0$. Thus $e_1(a) > 0$ for $0 < a \leq 5$. 
We claim that if $0 < a \leq 5$ and $\norm{(u/a,v)} \leq b$, then
\begin{equation}\label{eq:approx_by_K0}
    |K(u,v)| \leq |K_{\circ}(u,v)| + \e_0(a,u,v) ,
\end{equation}
where
\[
\e_0(a,u,v) := c(6, a) \frac{|n_0(u,v)|}{d_0(u,v)^{4}} 
\]
with 
\[
c(\alpha,a) := \frac{\alpha e_0(a)}{4 e_1(a)^{\frac{\alpha+4}{4}}}.
\]
Indeed, when $\norm{(u/a,v)} \leq b$, we have 
\begin{equation} \label{eq:d_est0}
\begin{split}
d(u,v)^4
&\geq  d_0(u,v)^4 - \Big(\frac{u^6}{4!} + \frac{u^2v^2}{2} + \frac{2 u^8}{8!} + \frac{2 v^4}{4!} \Big) \\
&\geq  d_0(u,v)^4 - \Big(\frac{u^6}{4!} + \frac{u^2v^2}{2} + \frac{2 u^6 (ab)^2}{8!} + \frac{2 |v|^3 b^2}{4!} \Big) \\
&\geq d_0(u,v)^4 \left(1 - e_0(a) d_0(u,v)^{2} \right),
\end{split}
\end{equation}
where we first used \eqref{eq:d_lower}, then used $|u| \leq ab$ and $|v| \leq b^2$ when $\norm{(u/a,v)} \leq b$, and finally used \eqref{eq:e0use}.
Furthermore, the mean value theorem gives
\[
\frac{1}{(1-t)^{\alpha/4}} \leq 1+ \frac{\alpha t}{4(1-t)^{\frac{\alpha+4}{4}}} \quad \text{if $\alpha > 0$ and $0 \leq t < 1$}.
\]
Applying this with $t = e_0(a) d_0(u,v)^{2} \leq e_0(a) d_0(ab,b^2)^{2} = 1-e_1(a) < 1$, and noting that 
\[
\frac{1}{(1-t)^\frac{\alpha+4}{4}} \leq \frac{1}{(1-e_0(a) d_0(ab,b^2)^{2})^{\frac{\alpha+4}{4}}} = \frac{1}{e_1(a)^{\frac{\alpha+4}{4}}},
\]
we get
\begin{equation} \label{eq:d_alpha_est}
\frac{1}{d(u,v)^{\alpha}} \leq \frac{1}{d_0(u,v)^{\alpha}} +\frac{c(\alpha,a)}{d_0(u,v)^{\alpha-2}}  \qquad \text{for $\alpha > 0$, $0 < a \leq 5$ and $\norm{(u/a,v)} \leq b$}.
\end{equation}
Inequality~\eqref{eq:approx_by_K0} now follows since $|n(u,v)| \leq |n_0(u,v)|$, which gives
\[
|K(u,v)| \leq \frac{|n_0(u,v)|}{d(u,v)^{6}},
\]
and since we can use \eqref{eq:d_alpha_est} with $\alpha = 6$ to bound $\frac{1}{d(u,v)^{6}}$.

We note also that 
\begin{equation} \label{eq:d_est2}
\frac{1}{d(u,v)} \leq \frac{1}{e_1(a)^{1/4} d_0(u,v)} \qquad \text{if $0 < a \leq 5$ and $\norm{(u/a,v)} \leq b$},
\end{equation}
because then we can apply \eqref{eq:d_est0} and use $1-e_0(a) d_0(u,v)^{2} \geq 1-e_0(a) d_0(ab,b^2)^{2} = e_1(a)$.

Using similar techniques, we can prove, for $0 < a \leq 5$ and $\norm{(u/a,v)} \leq b$,
\begin{equation} \label{eq:approx_by_Ku}
    |\frac{\partial K}{\partial u}(u,v)| \leq |\frac{\partial K_{\circ}}{\partial u}(u,v)| + \e_u(a,u,v)
\end{equation}
and
\begin{equation} \label{eq:approx_by_Kv}
    |\frac{\partial K}{\partial v}(u,v)| \leq |\frac{\partial K_{\circ}}{\partial v}(u,v)| + \e_v(a,u,v) ,
\end{equation}
where
\[
\e_u(a,u,v) := c(10, a) \frac{|n_{0u}(u,v)|}{d_0(u,v)^{8}} +  \frac{\frac{1451 |u^6 v| }{3072} + \frac{|u^2 v^3|}{2} +  \frac{v^4}{3} }{e_1(a)^{10/4} d_0(u,v)^{10}}, 
\]
and
\[
\e_v(a,u,v) := c(10, a) \frac{|n_{0v}(u,v)|}{d_0(u,v)^{8}} + \frac{\frac{683|u^7| }{15360}  + |u^3 v^2| + \frac{7|u v^3|}{6} }{e_1(a)^{10/4} d_0(u,v)^{10}}.
\]
This is because when $0 < a \leq 5$ and $\norm{(u/a,v)} \leq b$,
\[
\begin{split}
|\frac{\partial K}{\partial u}(u,v)| &\leq \frac{|n_{0u}(u,v)|}{d(u,v)^{10}}+ \frac{\frac{1451 |u^6 v| }{3072} + \frac{|u^2 v^3|}{2} +  \frac{v^4}{3} }{d(u,v)^{10}}  \\
&\leq  \frac{|n_{0u}(u,v)|}{d_0(u,v)^{10}}+ c(10,a) \frac{|n_{0u}(u,v)|}{d_0(u,v)^{8}}   +  \frac{\frac{1451 |u^6 v| }{3072} + \frac{|u^2 v^3|}{2} +  \frac{v^4}{3} }{e_1(a)^{10/4} d_0(u,v)^{10}},
\end{split}
\]
where in the first inequality we applied \eqref{eq:nu_est1}, and in second inequality we used \eqref{eq:d_alpha_est} with $\alpha = 10$ to estimate the first term, and used \eqref{eq:d_est2} to estimate the second term. This proves \eqref{eq:approx_by_Ku}. Similarly we have \eqref{eq:approx_by_Kv}. We also have, for $0 < a \leq 5$ and $\norm{(u/a,v)} \leq b$, that
\begin{equation} \label{eq:approx_by_Kuu}
    |\frac{\partial^2 K}{\partial u^2}(u,v)| \leq |\frac{\partial^2 K_{\circ}}{\partial u^2}(u,v)| + \e_{uu}(a,u,v)
\end{equation}
where
\[
\e_{uu}(a,u,v) := c(14,a) \frac{|n_{0uu}(u,v)|}{d_0(u,v)^{12}}  + \frac{\frac{1195 |u^9 v|}{2048} + \frac{|u^5 v^3|}{8} + \frac{5|u^3 v^4|}{2}  + \frac{15 |u v^5|}{4}}{ e_1(a)^{14/4} d_0(u,v)^{14}}.
\]
The numerical values of $e_0(a)$, $e_1(a)$ and $c(\alpha,a)$ can all be estimated using Mathematica once $\alpha$ and $a$ are given. Thus $\e_0$, $\e_u$, $\e_v$ and $\e_{uu}$ are completely concrete functions.

In the next two lemmas, we will prove two bounds for $\widehat{K}(m,n)$. The first one is good when $m$ is large; the second one is good when $n$ is large. Recall $b = \frac{2}{5}$.

\begin{Lem} \label{lem:Khat1}
Suppose $c_1 > 0$, $m$, $n$ are positive integers with the same parity and $m \geq c_1/b$. Then for any $0 < a \leq 5$,
\[
\begin{split}
|\widehat{K}(m,n)| \leq & \, c_1^3 I_0 \frac{n}{m^2} + 
 \frac{I_{1}}{c_1} + \frac{I_{2} + I_{3}}{c_1^2} \\
 & \quad + \frac{1}{m^2} \left( c_1^5 \frac{n}{m^2} \eps_{I0} + c_1 \eps_{I1} + \eps_{I2} - \frac{I_3}{b^2} + \eps_{I3} \log \frac{b m}{c_1} +  \eps_{I4}  \right)
\end{split}
\]
where $I_0, \dots, I_3$, $\eps_{I0},\dots,\eps_{I4}$ are constants depending on $a$, given by
\begin{alignat*}{3}
I_0 &:= \iint_{\norm{(u/a,v)} \leq 1} |u v K_{\circ}(u,v)| du dv, \qquad 
&&I_{1} := 2 \int_{-1}^1 |K_{\circ}(a,v)| dv, \\
I_{2} &:= 2 \int_{-1}^1 |\frac{\partial K_{\circ}}{\partial u} (a,v)| dv, 
&&I_{3} := \iint_{\norm{(u/a,v)} \geq 1} |\frac{\partial^2 K_{\circ}}{\partial u^2}(u,v)| du dv,
\end{alignat*}
and
\begin{align*}
\eps_{I0} &:= \iint_{\norm{(u/a,v)} \leq 1} |u v| \e_0(a,u,v) du dv, \\
\eps_{I1} &:= 2 \int_{-1}^1 \e_0(a,a,v) dv, \qquad
\eps_{I2} := 2 \int_{-1}^1 \e_u(a,a,v) dv, \\
\eps_{I3} &:= 2 a \int_{-1}^1 \e_{uu}(a,a,v) dv + 4 \int_{-a}^a \e_{uu}(a,u,1) du, \\
\eps_{I4} &:= \iint_{\substack{|u| \leq \pi, |v| \leq \pi/2 \\ \norm{(u/a,v)} > b}} |\frac{\partial^2 K}{\partial u^2}(u,v)| du dv.
\end{align*}
\end{Lem}

\begin{proof}
Define $R = c_1/m$. Our hypothesis guarantees that $R \leq b$. We split
\[
\begin{split}
-\widehat{K}(m,n) = &\, \iint_{\substack{|u| \leq \pi \\ |v| \leq \pi/2}} K(u,v) \sin(mu) \sin(nv) du dv \\
= &\, \iint_{\norm{(u/a,v)} \leq R} + \iint_{\substack{|u| \leq \pi, |v| \leq \pi/2 \\ \norm{(u/a,v)} > R}} K(u,v) \sin(mu) \sin(nv) du dv.
\end{split}
\]
We estimate the first term by putting absolute value inside, and integrate by parts in the second term using $\sin(mu)=-\frac{1}{m} \frac{d}{du} \cos(mu)$. Note that the set $\{(u,v) \colon \norm{(u/a,v)} = R\}$ is the boundary of a rectangle. Thus
\[
\begin{split}
& \iint_{\substack{|u| \leq \pi, |v| \leq \pi/2 \\ \norm{(u/a,v)} > R}} K(u,v) \sin(mu) \sin(nv) du dv \\
& = \, -\frac{2}{m} \int_{-R^2}^{R^2} K(aR,v) \cos(m a R) \sin(nv) dv  + \frac{1}{m} \iint_{\substack{|u| \leq \pi, |v| \leq \pi/2 \\ \norm{(u/a,v)} > R}} \frac{\partial K}{\partial u}(u,v) \cos(mu) \sin(nv) du dv.
\end{split}
\]
We integrate by parts once more for the last double integral, using $\cos(mu)=\frac{1}{m} \frac{d}{du} \sin(mu)$, and obtain
\[
\begin{split}
& \frac{1}{m} \iint_{\substack{|u| \leq \pi, |v| \leq \pi/2 \\ \norm{(u/a,v)} > R}} \frac{\partial K}{\partial u}(u,v) \cos(mu) \sin(nv) du dv \\
& = \, \frac{2}{m^2} \int_{-R^2}^{R^2} \frac{\partial K}{\partial u}(aR,v) \sin(m a R) \sin(nv) dv  - \frac{1}{m^2} \iint_{\substack{|u| \leq \pi, |v| \leq \pi/2 \\ \norm{(u/a,v)} > R}} \frac{\partial^2 K}{\partial u^2}(u,v) \sin(mu) \sin(nv) du dv.
\end{split}
\]
Taking absolute values, and using $|\sin(mu) \sin(nv)| \leq m n |uv|$, we get
\[
\begin{split}
|\widehat{K}(m,n)| & \leq \, m n \iint_{\norm{(u/a,v)} \leq R} |u v K(u,v)| du dv + \frac{2}{m} \int_{-R^2}^{R^2} |K(aR,v)| dv + \frac{2}{m^2} \int_{-R^2}^{R^2} |\frac{\partial K}{\partial u}(aR,v)| dv  \\
& \quad + \frac{1}{m^2} \iint_{R < \norm{(u/a,v)} \leq b} |\frac{\partial^2 K}{\partial u^2}(u,v)| du dv + \frac{1}{m^2} \iint_{\substack{|u| \leq \pi, |v| \leq \pi/2 \\ \norm{(u/a,v)} > b}} |\frac{\partial^2 K}{\partial u^2}(u,v)| du dv.
\end{split}
\]
The first four integrals are all contained in $\{(u,v) \colon \norm{(u/a,v)} \leq b\}$, so our earlier estimates \eqref{eq:approx_by_K0}, \eqref{eq:approx_by_Ku} and \eqref{eq:approx_by_Kuu} apply. They allow us to bound $K$ and its derivatives by those of $K_{\circ}$, up to some errors that we control. 
Then we obtain
\[
\begin{split}
|\widehat{K}(m,n)| & \leq \, m n \iint_{\norm{(u/a,v)} \leq R} |u v K_{\circ}(u,v)| du dv + \frac{2}{m} \int_{-R^2}^{R^2} |K_{\circ}(aR,v)| dv + \frac{2}{m^2} \int_{-R^2}^{R^2} |\frac{\partial K_{\circ}}{\partial u}(aR,v)| dv  \\
& \quad + \frac{1}{m^2} \iint_{R < \norm{(u/a,v)} \leq b} |\frac{\partial^2 K_{\circ}}{\partial u^2}(u,v)| du dv \\
& \quad +m n \iint_{\norm{(u/a,v)} \leq R} |u v| \e_0(a,u,v) du dv + \frac{2}{m} \int_{-R^2}^{R^2} \e_0(a,aR,v) dv + \frac{2}{m^2} \int_{-R^2}^{R^2} \e_u(a,aR,v) dv  \\
& \quad + \frac{1}{m^2} \iint_{R < \norm{(u/a,v)} \leq b} \e_{uu}(a,u,v) du dv + \frac{\eps_{I4}}{m^2}.
\end{split}
\]
All the integrals on the right hand side above are now integrals of homogeneous functions of $(u,v)$. Since $|u v K_{\circ}(u,v)|$ is homogeneous of degree $0$, a change of variables $(u,v) \mapsto (R u, R^2 v)$ give
\[
 m n \iint_{\norm{(u/a,v)} \leq R} |u v K_{\circ}(u,v)| du dv =  m n R^3 \iint_{\norm{(u/a,v)} \leq 1} |u v K_{\circ}(u,v)| du dv = m n R^3 I_0.
\]
Similarly, the integrals involving $K_{\circ}$, $\frac{\partial K_{\circ}}{\partial u}$, $\e_0$ and $\e_u$ are homogeneous functions of $R$. Next, we rewrite the integral involving $\frac{\partial^2 K_{\circ}}{\partial u^2}$ by a change of variables:
\[
\begin{split}
&\quad \frac{1}{m^2} \iint_{R < \norm{(u/a,v)} \leq b} |\frac{\partial^2 K_{\circ}}{\partial u^2}(u,v)| du dv  \\
&= \frac{1}{m^2} \left( \iint_{\norm{(u/a,v)} > R} |\frac{\partial^2 K_{\circ}}{\partial u^2}(u,v)| du dv - \iint_{\norm{(u/a,v)} > b} |\frac{\partial^2 K_{\circ}}{\partial u^2}(u,v)| du dv \right) \\
&= I_3 \left(\frac{1}{m^2 R^2} - \frac{1}{m^2 b^2} \right).
\end{split}
\]
Finally, 
by the Fundamental Theorem of Calculus and the chain rule,
\[
\frac{d}{dt} \iint_{1 < \norm{(u/a,v)} \leq t} \e_{uu}(a,u,v) du dv = \, 2 a \int_{-t^2}^{t^2} \e_{uu}(a,at,v) dv + 4 t \int_{-at}^{at} \e_{uu}(a,u,t^2) du,
\]
which, by homogeneity, is equal to
\[
\frac{1}{t} \Big( 2 a \int_{-1}^{1} \e_{uu}(a,a,v) dv + 4 \int_{-a}^{a} \e_{uu}(a,u,1) du \Big) = \frac{\eps_{I3}}{t}.
\]
Thus
\[
\begin{split}
\iint_{R < \norm{(u/a,v)} \leq b} \e_{uu}(a,u,v) du dv
= &\, \iint_{1 < \norm{(u/a,v)} \leq b/R} \e_{uu}(a,u,v) du dv \\
= &\, \int_1^{b/R} \frac{d}{dt} \iint_{1 < \norm{(u/a,v)} \leq t} \e_{uu}(a,u,v) du dv dt \\
= &\, \int_1^{b/R} \frac{\eps_{I3}}{t} dt = \eps_{I3} \log(\frac{b}{R}),
\end{split}
\]
and 
putting all these together, one sees that 
\[
\begin{split}
|\widehat{K}(m,n)| 
\leq & \, m n R^3 I_0 + \frac{I_{1}}{m R} + \frac{I_{2}+I_{3}}{m^2 R^2} \\ & \quad +m n R^5 \eps_{I0} + \frac{\eps_{I1}}{m} R + \frac{\eps_{I2}}{m^2} - \frac{I_3}{m^2 b^2} + \frac{\eps_{I3} \log(b/R)}{m^2} + \frac{\eps_{I4}}{m^2}.
\end{split}
\]
Remembering $R = c_1/m$ gives the estimate in Lemma \ref{lem:Khat1}.
\end{proof}

\begin{Lem}  \label{lem:Khat2}
Suppose $c_2 > 0$, $m$, $n$ are positive integers with the same parity and $n \geq (c_2/b)^2$. Then for any $0 < a \leq 5$, 
\[
\begin{split}
|\widehat{K}(m,n)| \leq &\,  c_2^3 J_0 \frac{m}{\sqrt{n}} + 
 \frac{J_{1} +J_{2}}{c_2^2} \\
 & \quad + \frac{1}{n} \left( c_2^5  \eps_{J0} \frac{m}{\sqrt{n}}  +  \eps_{J1} - \frac{J_2}{b^2} + \eps_{J2}  \log \frac{b \sqrt{n}}{c_2} + \eps_{J3} \right)
 \end{split}
\]
where $J_0,\dots,J_2$ and $\eps_{J0},\dots,\eps_{J3}$ are constants depending on $a$, given by
\begin{align*}
J_{0} &:= \iint_{\norm{(u/a,v)} \leq 1} |u v K_{\circ}(u,v)| du dv, \\
J_{1} &:= 2 \int_{-a}^a |K_{\circ}(u,1)| du, \qquad 
J_{2} := \iint_{\norm{(u/a,v)} \geq 1} |\frac{\partial K_{\circ}}{\partial v}(u,v)| du dv,
\end{align*}
and
\begin{align*}
\eps_{J0} &:= \iint_{\norm{(u/a,v)} \leq 1} |u v | \e_0(a,u,v) du dv, \\
\eps_{J1} &:= 2 \int_{-a}^a \e_0(a,u,1) du \\
\eps_{J2} &:= 2 a \int_{-1}^1 \e_v(a,a,v) dv + 4 \int_{-a}^a \e_v(a,u,1) du, \\
\eps_{J3} &:= \iint_{\substack{|u| \leq \pi, |v| \leq \pi/2 \\ \norm{(u/a,v)} > b}} |\frac{\partial K}{\partial v}(u,v)| du dv.
\end{align*}
\end{Lem}

\begin{proof}
The proof is similar to the previous lemma (and slightly easier); one only needs to split the integral into two parts, one where $\norm{(u/a,v)} \leq R$, and another where $\norm{(u/a,v)} > R$, with $R := c_2/\sqrt{n} < b$; one then integrates by parts once in $v$ for the second integral, using $\sin(n v) = -\frac{1}{n} \frac{d}{dv} \cos(n v)$. 

More precisely, define $R = c_2/\sqrt{n}$. Our hypothesis guarantees that $R \leq b$. We split\[
\begin{split}
-\widehat{K}(m,n) = &\, \iint_{\substack{|u| \leq \pi \\ |v| \leq \pi/2}} K(u,v) \sin(mu) \sin(nv) du dv \\
= &\, \iint_{\norm{(u/a,v)} \leq R} + \iint_{\substack{|u| \leq \pi, |v| \leq \pi/2 \\ \norm{(u/a,v)} > R}} K(u,v) \sin(mu) \sin(nv) du dv.
\end{split}
\]
We estimate the first term by putting absolute value inside, and integrate by parts in the second term using $\sin(n v) = -\frac{1}{n} \frac{d}{dv} \cos(n v)$. Note that the set $\{(u,v) \colon \norm{(u/a,v)} = R\}$ is the boundary of a rectangle. Thus
\[
\begin{split}
& \iint_{\substack{|u| \leq \pi, |v| \leq \pi/2 \\ \norm{(u/a,v)} > R}} K(u,v) \sin(mu) \sin(nv) du dv \\
& = -\, \frac{2}{n} \int_{-aR}^{aR} K(u,R^2) \sin(m u) \cos(n R^2) du  + \frac{1}{n} \iint_{\substack{|u| \leq \pi, |v| \leq \pi/2 \\ \norm{(u/a,v)} > R}} \frac{\partial K}{\partial v}(u,v) \sin(mu) \cos(nv) du dv.
\end{split}
\]
Thus
\[
\begin{split}
|\widehat{K}(m,n)| & \leq \, m n \iint_{\norm{(u/a,v)} \leq R} |u v K(u,v)| du dv + \frac{2}{n} \int_{-aR}^{aR} |K(u,R^2)| du \\
& \quad + \frac{1}{n} \iint_{R < \norm{(u/a,v)} \leq b} |\frac{\partial K}{\partial v}(u,v)| du dv + \frac{1}{n} \iint_{\substack{|u| \leq \pi, |v| \leq \pi/2 \\ \norm{(u/a,v)} > b}} |\frac{\partial K}{\partial v}(u,v)| du dv.
\end{split}
\]

The first three integrals are all contained in $\{(u,v) \colon \norm{(u/a,v)} \leq b\}$, so our earlier estimates \eqref{eq:approx_by_K0} and \eqref{eq:approx_by_Kv}  apply.  Thus
\[
\begin{split}
|\widehat{K}(m,n)| & \leq \, m n \iint_{\norm{(u/a,v)} \leq R} |u v K_{\circ}(u,v)| du dv + \frac{2}{n} \int_{-aR}^{aR} |K_{\circ}(u,R^2)| du + \frac{1}{n} \iint_{R < \norm{(u/a,v) \leq b} } |\frac{\partial K_{\circ}}{\partial v}(u,v)| du dv \\
& \quad +m n \iint_{\norm{(u/a,v)} \leq R} |u v| \e_0(a,u,v) du dv + \frac{2}{n} \int_{-aR}^{aR} \e_0(a,u,R^2) dv \\
& \quad + \frac{1}{n} \iint_{R < \norm{(u/a,v)} \leq b} \e_{v}(a,u,v) du dv  + \frac{\eps_{J3}}{m^2}.
\end{split}
\]
All the integrals on the right hand side above are now integrals of homogeneous functions of $(u,v)$. Thus the integrals are homogeneous functions of $R$, with the exception of the one involving $\e_{v}$ because $\e_{v}$ is homogeneous of degree $-3$. But then, by the Fundamental Theorem of Calculus and the chain rule,
\[
\frac{d}{dt} \iint_{1 < \norm{(u/a,v)} \leq t} \e_{v}(a,u,v) du dv = \, 2 a \int_{-t^2}^{t^2} \e_{v}(a,at,v) dv + 4 t \int_{-at}^{at} \e_{v}(a,u,t^2) du,
\]
which, by homogeneity, is equal to
\[
\frac{1}{t} \Big( 2 a \int_{-1}^{1} \e_{v}(a,a,v) dv + 4 \int_{-a}^{a} \e_{v}(a,u,1) du \Big) = \frac{\eps_{J2}}{t}.
\]
Thus
\[
\begin{split}
\iint_{R < \norm{(u/a,v)} \leq b} \e_{v}(a,u,v) du dv
= &\, \iint_{1 < \norm{(u/a,v)} \leq b/R} \e_{v}(a,u,v) du dv \\
= &\, \int_1^{b/R} \frac{d}{dt} \iint_{1 < \norm{(u/a,v)} \leq t} \e_{v}(a,u,v) du dv dt \\
= &\, \int_1^{b/R} \frac{\eps_{J2}}{t} dt = \eps_{J2} \log(\frac{b}{R}),
\end{split}
\]
and by scaling all other integrals similarly one sees that 
\[
\begin{split}
|\widehat{K}(m,n)| 
\leq & \, m n R^3 J_0 + \frac{J_{1}+J_{2}}{n R^2} \\ & \quad +m n R^5 \eps_{J0} + \frac{\eps_{J1}}{n} - \frac{J_2}{n b^2} + \frac{\eps_{J2} \log(b/R)}{n} + \frac{\eps_{J3}}{n}.
\end{split}
\]
Remembering $R = c_2/\sqrt{n}$ gives the estimate in Lemma \ref{lem:Khat2}.
\end{proof}

We apply the above two lemmas with $a = 2$, $c_1:=\frac{8}{5}$ and $a = 1$, $c_2:=\frac{4\sqrt{2}}{5}$ respectively:

\begin{Cor} \label{cor:khat1}
Suppose $m,n$ are positive integers with the same parity. If $m \geq 4$, then
\begin{equation} \label{eq:khatest1}
|\widehat{K}(m,n)| \leq \frac{7168}{625} \frac{n}{m^2} + 5 + \frac{1}{m^2} \Big(\frac{783}{40} + \frac{98304}{3125} \frac{n}{m^2} + 192 \log(\frac{m}{4})\Big).
\end{equation}
\end{Cor}

\begin{proof}
We apply Lemma~\ref{lem:Khat1} with $a = 2$, so that Mathematica returns  expressions for exact values for $I_0, I_1, I_2, I_3, \eps_{I0}, \eps_{I1}, \eps_{I2}, \eps_{I3}$ and a numerical estimate for $\eps_{I4}$ (see Appendix A
). We then have
\[
\begin{split}
I_0 &= \frac{32}{3}  \Big( \log(\frac{1+\sqrt{5}}{2}) - \frac{1}{32} \sqrt{\frac{2}{\pi}}  \Gamma(-\frac{1}{4})^2 + \,_2F_1(-\frac{1}{4},\frac{1}{2},\frac{3}{4},-\frac{1}{4}) -\, _2F_1(\frac{1}{2},\frac{3}{4},\frac{7}{4},-4)\Big) \\
I_1 &= 4 - \frac{8}{\sqrt{5}}, \qquad
I_2 = 2 - \frac{12}{5 \sqrt{5}}, \qquad 
I_3 = 6 + 4 \sqrt{2} - \frac{12}{5 \sqrt{5}}, 
\end{split}
\]
so
\[
I_0 < \frac{14}{5}, \quad
I_1 < \frac{1}{2}, \quad 
I_2 < 1, \quad 
\frac{21}{2} < I_3 < 11.
\]
Next, recall $b = 2/5$, so
\[
\begin{split}
e_0(2) &= \frac{1}{75} \sqrt{\frac{68935125533 + 1312025 \sqrt{1102101}}{55105058}}, \\
e_1(2) &= 1 - \frac{4}{5\sqrt{5}} e_0(2),
\end{split}
\]
\[
\begin{split}
\eps_{I0} &= \frac{3}{5} \frac{e_0(2)}{e_1(2)^{5/2}} \Big(16 + \pi + 4 \arctan(3) - 32 \arctan (\frac{1}{2}) - 2 \log 5 \Big), \\
\eps_{I1} &= \frac{6 e_0(2)}{e_1(2)^{5/2}} \log \frac{5}{4}, \\
\eps_{I2} &= (6 + 5 \log \frac{4}{5}) \frac{e_0(2)}{e_1(2)^{7/2}} +  \frac{1}{e_1(2)^{5/2} } \Big(\frac{2219}{288} - \frac{307}{20 \sqrt{5}} + \frac{4}{3} \log \frac{1+\sqrt{5}}{2} \Big),  \\
\eps_{I3} &= 105 \frac{e_0(2)}{e_1(2)^{9/2}} + \frac{9473}{240 \, e_1(2)^{7/2}}.
\end{split}
\]
Together with the numerical value of $\eps_{I4}$ we get
\[
\eps_{I0} < 3, \quad \eps_{I1} < 2, \quad \eps_{I2} < 7, \quad \eps_{I3} < 192, \quad \eps_{I4} < 75.
\]
Finally we set $c_1 = \frac{8}{5}$ so that $c_1/b = 4 \leq m$. 
\end{proof}

\begin{Cor} \label{cor:khat2}
Suppose $m,n$ are positive integers with the same parity. If $n \geq 8$, then
\begin{equation} \label{eq:khatest2}
|\widehat{K}(m,n)| \leq \frac{4832 \sqrt{2}}{3125} \frac{m}{\sqrt{n}} + \frac{245}{32} + \frac{1}{n} \Big( \frac{53}{8} + \frac{16384 \sqrt{2}}{15625} \frac{m}{\sqrt{n}} + 19 \log(\frac{n}{8}) \Big).
\end{equation}
\end{Cor}


\begin{proof}
We apply Lemma~\ref{lem:Khat2} with $a = 1$, so that Mathematica returns expressions for exact values for $J_0, J_1, J_2, \eps_{J0}, \eps_{J1}, \eps_{J2}$ and a numerical estimate for $\eps_{J3}$ (see Appendix A).
Thus we have:
\[
\begin{split}
J_0 &= \frac{1}{3} \Big(4 \log(2+\sqrt{5}) - \sqrt{ \frac{2}{\pi} } \Gamma(-\frac{1}{4})^2 + 16 \, _2F_1(-\frac{1}{4},\frac{1}{2},\frac{3}{4},-4) - 4 \,_2F_1(\frac{1}{2},\frac{3}{4},\frac{7}{4},-\frac{1}{4}) \Big), \\
J_1 &= \frac{4}{\sqrt{5}}, \qquad
J_2 = \frac{32}{3\sqrt{3}} + \frac{4}{\sqrt{5}},
\end{split}
\]
so
\[
J_0 < \frac{151}{100}, \quad J_1 < \frac{9}{5}, \quad \frac{79}{10} < J_2 < 8.
\]
Next, recall $b = 2/5$, so
\[
\begin{split}
e_0(1) &= \frac{1}{75} \sqrt{\frac{137909641801 + 2625550 \sqrt{1102731}}{110273101}} \\
e_1(1) &= 1-\frac{2}{5 \sqrt{5}} e_0(1)
\end{split}
\]
\[
\begin{split}
\eps_{J0} &= \frac{3}{5} \frac{e_0(1)}{e_1(1)^{5/2}} \Big( 2 + 3 \arctan(2) - 2 \log 5 \Big) \\
\eps_{J1} &= 6 \arctan(\frac{1}{2})  \frac{e_0(1)}{e_1(1)^{5/2}}  \\
\eps_{J2} &= 5 \Big(4 \sqrt{2} + \pi - 4 \arctan(\frac{1}{\sqrt{2}}) \Big)  \frac{e_0(1)}{e_1(1)^{7/2}} + \frac{4843}{360 \, e_1(1)^{5/2}}.
\end{split}
\]
Together with the numerical value of $\eps_{J3}$, we get
\[
\eps_{J0} < \frac{4}{5}, \quad \eps_{J1} < 2, \quad \eps_{J2} < 38, \quad \eps_{J3} < 54.
\]
Finally we set $c_2 = \frac{4\sqrt{2}}{5}$ so that $(c_2/b)^2 = 8 \leq n$. 
\end{proof}
\begin{rem}The Mathematica codes used in Corollaries \ref{cor:khat1} and \ref{cor:khat2} are enclosed in Appendix A for interested readers. 
\end{rem}

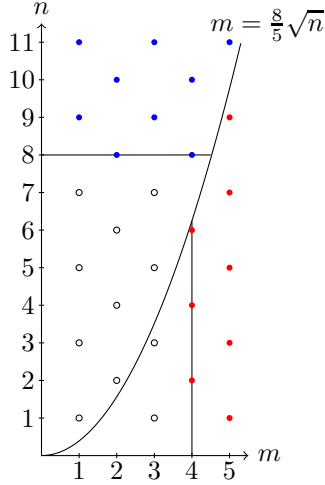
\begin{figure}
\begin{tikzpicture}[domain=0:5.5, x = 0.5cm, y = 0.5cm]
\draw[->] (0,0) -- (5.5,0);
\draw[->] (0,0) -- (0,11.5);
\draw (5.5,0) node [right] {$m$};
\draw (0,11.5) node [above] {$n$};
\draw (6,10.8) node [above] {$m = \frac{8}{5} \sqrt{n}$};
\foreach \x in {1,...,5} {%
    \draw ($(\x,0) + (0,-1pt)$) -- ($(\x,0) + (0,1pt)$)
        node [below] {$\x$};
}
\foreach \y in {1,...,11} {%
    \draw ($(0,\y) + (-1pt,0)$) -- ($(0,\y) + (1pt,0)$)
        node [left] {$\y$};
}
\draw[black]   plot[smooth,domain=0:5.3] (\x, {\x^2/(8/5)^2});
\draw[black] (0,8) -- (4.5,8);
\draw[black] (4,0)--(4,6.25);
\foreach \i in {1,3,...,9}
\fill[red] (5,\i) circle (0.04cm);
\foreach \i in {2,4,6}
\fill[red] (4,\i) circle (0.04cm);
\foreach \i in {1,3,5}
\fill[blue] (\i,11) circle (0.04cm);
\foreach \i in {2,4}
\fill[blue] (\i,8) circle (0.04cm);
\foreach \i in {1,3}
\fill[blue] (\i,9) circle (0.04cm);
\foreach \i in {2,4}
\fill[blue] (\i,10) circle (0.04cm);
\foreach \i in {1,3,...,7}
\draw (1,\i) circle (0.04cm);
\foreach \i in {2,4,6}
\draw (2,\i) circle (0.04cm);
\foreach \i in {1,3,5,7}
\draw (3,\i) circle (0.04cm);
\end{tikzpicture}
\caption{The points for which $|\widehat{K}(m,n)|$ can be bounded by Corollaries~\ref{cor:khat1} and~\ref{cor:khat2} are coloured in red and blue respectively.
} 
\label{fig3}
\end{figure}

From Corollary~\ref{cor:khat1}, if $m \geq 4$ and $\frac{8}{5} \sqrt{n} \leq m$, then 
\begin{equation} \label{eq:used_m_over_rootn_bdd1}
|\widehat{K}(m,n)| \leq  \frac{237}{25} + \frac{1}{m^2} \Big( \frac{31863}{1000} + 192 \log(\frac{m}{4})\Big).
\end{equation}
The right hand side is an increasing function of $m$ on $[4,5]$ and a decreasing function of $m$ on $[6,\infty)$, and is $<4 \pi$ at both $m = 5$ and $m = 6$. Thus we have $|\widehat{K}(m,n)| < 4\pi$ if $m \geq 4$ and $m \geq \frac{8}{5} \sqrt{n}$. This shows $|\widehat{K}(m,n)| < 4\pi$ for $(m,n)$ coloured in red in Figure~\ref{fig3}. 

From Corollary~\ref{cor:khat2}, if $n \geq 8$ and $m \leq \frac{8}{5} \sqrt{n}$, then 
\begin{equation} \label{eq:used_m_over_rootn_bdd2}
|\widehat{K}(m,n)| \leq \frac{38656\sqrt{2}}{15625} + \frac{245}{32} +\frac{1}{n} \Big( \frac{131072 \sqrt{2}}{78125} + \frac{53}{8} + 19 \log(\frac{n}{8}) \Big)
\end{equation}
The right hand side is an increasing function of $n$ on $[8,13]$ and decreasing on $[14,\infty)$. Since this function of $n$ is $<4 \pi$ at both $n=13$ and $n=14$, we have $|\widehat{K}(m,n)| < 4\pi$ if $n \geq 8$ and $m \leq \frac{8}{5} \sqrt{n}$. Such pairs $(m,n)$ are pictured in blue in Figure~\ref{fig3}.


It remains to verify that $|\widehat{K}(m,n)| < 4\pi$ for the 11 points in Figure~\ref{fig3} that are not coloured, i.e. when $1 \leq m \leq 3$, $1 \leq n \leq 7$ and $(m,n)$ are of the same parity. We achieved this numerically using Mathematica (in fact, $|\widehat{K}(m,n)| < 5.8$ for all $(m,n)$ that is not coloured). For interested readers, see Appendix B for the Mathematica codes.

\section{The weak solution to the Neumann problem} \label{sect5}

In this section we prove Theorem~\ref{thm:NeumannL2}.

\subsection{Functional Analysis}
Recall that the Folland--Stein space $S^{1,2}(\Omega)$ is the space of all functions $u \in L^2(\Omega)$ whose weak derivatives $\nablab u$ are in $L^2(\Omega)$. 
It is known~\cite[Theorem~1.3]{Nhieu2001} that $C^\infty(\overline{\Omega}) \subset S^{1,2}(\Omega)$ is a dense subspace and that the restriction map
\begin{equation*}
 C^\infty(\overline{\Omega}) \ni u \mapsto u\rvert_{\Sigma} \in C^\infty(\Sigma)
\end{equation*}
is continuous.
Therefore the trace operator $\Tr \colon S^{1,2}(\Omega) \to L^2(\Sigma)$ is defined.
It is known~\cite[Theorem~1.4]{Nhieu2001} that $\Tr\bigl(S^{1,2}(\overline{\Omega})\bigr) \subset L^3(\Sigma)$.

Since $\Sigma$ is $p$-minimal, we see that
\begin{equation*}
 \mathcal{F}(u) := \int_\Omega \left( \lvert \nablab u \rvert^2 + \frac{R}{4}\lvert u \rvert^2 \right) \, \theta \wedge d\theta
\end{equation*}
is the continuous extension of $\mathcal{F} \colon C^\infty(\Omega) \to \mathbb{R}$ to $S^{1,2}(\Omega)$.
Moreover, there is a minimizer $u_0 \in S^{1,2}(\Omega)$ of
\begin{equation*}
 \mu_{1} := \left\{ \mathcal{F}(u) \mathrel{}:\mathrel{} u \in S^{1,2}(\Omega) , \oint_{\Sigma} \lvert \Tr(u) \rvert^2 \, d\sigma = 1 \right\} .
\end{equation*}
Necessarily $u_0$ is a weak solution of
\begin{equation*}
 \begin{cases}
  Lu = 0 , & \text{in $\Omega$} , \\
  \nabla_\nu u = -\mu_{1}\Tr u , & \text{on $\Sigma$};
 \end{cases}
\end{equation*}
i.e.
\begin{equation*}
 0 = \int_\Omega \left( \langle \nablab u_0 , \nablab \phi \rangle + \frac{1}{2}u_0\phi \right) \, \theta \wedge d\theta + \oint_\Sigma \mu_{1}\Tr(u_0)\phi \, d\sigma
\end{equation*}
for all $\phi \in C^\infty(\overline{\Omega})$.
In this section we show that $u_0 \in C^\infty(\overline{\Omega})$.

\subsection{Potential theory}
Let $f \in L^2(\Sigma,d\sigma)$. Recall the single layer potential 
\[
\mathcal{S}f(\zeta) =  \int_{\eta \in \Sigma} f(\eta) G(\zeta,\eta) d\sigma(\eta), \quad \zeta \in \overline{\Omega}.
\]
We have $\mathcal{S}f \in C^{\infty}(\Omega)$. By \eqref{eq:Sf_uv_nice}, we also have
\[
\mathcal{S}f(\zeta) = \frac{\sqrt{2}}{8\pi} \iint_{\substack{|u| \leq \pi \\ |v| \leq \pi/2}} F(u'-u,v'-v) k(\frac{t}{\sqrt{2}},u,v) du dv
\]
if $\zeta = \Phi(t,u',v')$ and $F(u,v) := f(\Phi(0,u,v))$.

\begin{prop}
If $f \in L^2(\Sigma,d\sigma)$ then  $\mathcal{S}f \in S^{1,2}(\Omega)$.
\end{prop}

\begin{proof}
 Note that 
\[
\sup_{t \in [0,\pi/100]} \iint_{\substack{|u| \leq \pi \\ |v| \leq \pi/2}} |k(\frac{t}{\sqrt{2}},u,v)| du dv \lesssim 1
\]
by the bound for $k$ in Lemma~\ref{lem:Kt}.
Using Young's convolution inequality on $\R^2 / \Lambda$, we see that
\[
\sup_{t \in [0,\pi/100]} \|\mathcal{S}f\|_{L^2(\Sigma_t)} \lesssim \|f\|_{L^2(\Sigma)}.
\]
Thus $\mathcal{S}f \in L^2(\Omega)$. Furthermore, for $\zeta \in \Omega$, we have
\[
Z_{\bar{1}} \mathcal{S}f(\zeta) = \int_{\eta \in \Sigma} f(\eta) Z_{\bar{1}} G(\zeta,\eta) d\sigma(\eta).
\]
But
\[
Z_{\bar{1}} G(\zeta,\eta) = \frac{1}{16 \pi} \frac{(\zeta^2 \eta^1 - \zeta^1 \eta^2)(1 - \zeta \cdot \overline{\eta})}{|1-\zeta \cdot \overline{\eta}|^3}.
\]
If $\zeta = \Phi(t,u',v')$ and $\eta = \Phi(0,u,v)$, then
\[
\begin{split}
\zeta^2 \eta^1 - \zeta^1 \eta^2 
= & \cos(\frac{\pi}{4}) \sin(\frac{\pi}{4}+\frac{t}{\sqrt{2}}) e^{i(v'-u'+v+u)} - \sin(\frac{\pi}{4})\cos(\frac{\pi}{4}+\frac{t}{\sqrt{2}}) e^{i(v'+u'+v-u)} \\
= & e^{i(v'+v)}  \Big[ \sin (\frac{t}{\sqrt{2}}) \cos(u'-u) - i \cos(\frac{t}{\sqrt{2}}) \sin(u'-u) \Big] ,
\end{split}
\]
so
\[
\begin{split}
Z_{\bar{1}} G(\zeta,\eta) & = \, \frac{e^{2iv'}}{16 \pi} \Big[ \sin (\frac{t}{\sqrt{2}}) \cos(u'-u) - i \cos(\frac{t}{\sqrt{2}}) \sin(u'-u) \Big]  \\
& \qquad \times \frac{\Big[ \Big(\cos(v'-v)  - \cos(\frac{t}{\sqrt{2}}) \cos(u'-u) \Big) -i \Big(\sin(v'-v) - \sin(\frac{t}{\sqrt{2}}) \sin(u'-u) \Big) \Big]}{g(t,u'-u,v'-v)^{3/2}} \\
& = \, -\frac{e^{2iv'}}{16 \pi} i \cos(\frac{t}{\sqrt{2}}) \frac{\sin(u'-u) \Big[ \Big(\cos(v'-v)  - \cos(\frac{t}{\sqrt{2}}) \cos(u'-u) \Big) -i \sin(v'-v) \Big]}{g(t,u'-u,v'-v)^{3/2}} \\
&\quad + O\Big(\frac{t}{g(t,u'-u,v'-v)}\Big)
\end{split}
\]
The numerator of the first term is even in $t$, but the denominator is not. Thus we symmetrize the first term by making it even in $t$, which gives an error that is also in $O\Big(\frac{t}{g(t,u'-u,v'-v)}\Big)$. In fact, the error term $\mathcal E$ is given by
\[\begin{split}
\mathcal E & = \, C
 \sin (u'-u) [ (\cos (v'-v) - \cos t \cos (u'-u) ) - i \sin (v'-v) ] \\
 & \quad \times \Big( \frac{1}{g(t,u'-u,v'-v)^{3/2}} - \frac{1}{g(-t,u'-u,v'-v)^{3/2}} \Big)\\
 & \leq C
 |u'-u| [ (\cos (v'-v) - \cos t \cos (u'-u) )^2 +  \sin ^2(v'-v)]^{1/2} \\
 & \quad \times \frac{1}{g(t^*,u'-u,v'-v)^{5/2}}
 \frac{\partial g}{\partial t} (t^*,u'-u,v'-v)\cdot t.
 \end{split}
\]
Plugging  
\[
\frac{\partial g}{\partial t} (t^*,u'-u,v'-v)
=g^{1/2}(t^*,u'-u,v'-v) \cdot O(|u'-u|+ |t|),
\]
where $t^*\in [-t, t]$,
into the error term, we obtain 
\[
\mathcal E= O\Big(\frac{t}{g(t,u'-u,v'-v)}\Big).\]

To summarize, if $\zeta = \Phi(t,u',v')$ and $\eta = \Phi(0,u,v)$, we get
\[
Z_{\bar{1}} G(\zeta,\eta) = -\frac{e^{2iv'}}{32 \pi} i \cos(\frac{t}{\sqrt{2}}) K_{t/\sqrt{2}}(u'-u,v'-v) + O\Big(\frac{t}{g(t,u'-u,v'-v)}\Big)
\]
where
\[
K_t(u,v) := \sin u [ (\cos v - \cos t \cos u) - i \sin v ] \Big( \frac{1}{g(t,u,v)^{3/2}} + \frac{1}{g(-t,u,v)^{3/2}} \Big).
\]
It follows that
\[
\begin{split}
Z_{\bar{1}} \mathcal{S} f(\Phi(t,u',v')) 
= & -\frac{e^{2iv'}}{32\pi} i \cos(\frac{t}{\sqrt{2}}) \iint_{\substack{|u| \leq \pi \\ |v| \leq \pi/2}} F(u'-u,v'-v) K_{t/\sqrt{2}}(u,v)  du dv \\
& \quad + O \left( \iint_{\substack{|u| \leq \pi \\ |v| \leq \pi/2}} |F(u'-u,v'-v)| \frac{t}{g(t,u,v)} du dv\right)
\end{split}
\]
is essentially the sum of the convolutions of $F(u,v) := f(\Phi(0,u,v))$ with two kernels on $\R^2 / \Lambda$.
The second term can be handled using
\[
\left\| \iint_{\substack{|u| \leq \pi \\ |v| \leq \pi/2}} F(u'-u,v'-v) \frac{t}{g(t,u,v)} du dv \right\|_{L^2(du' dv')} \lesssim \|F\|_{L^2} ,
\]
where we applied Young's convolution inequality on $\R^2 / \Lambda$ and the bound
\[
\begin{split}
    \iint_{\substack{|u| \leq \pi \\ |v| \leq \pi/2}} \frac{t}{g(t,u,v)} du dv
&\lesssim \iint_{\substack{|u| \leq \pi \\ |v| \leq \pi/2}} \frac{t}{(|t|+|u|+|v|^{1/2})^4} du dv \\
&\leq \iint_{\R^2} \frac{1}{(1+|u|+|v|^{1/2})^4} du dv < \infty
\end{split}
\]
uniformly in $t \in (0,\pi/100]$. We claim that 
\[
\left\| \iint_{\substack{|u| \leq \pi \\ |v| \leq \pi/2}} F(u'-u,v'-v) K_t(u,v) du dv\right\|_{L^2(du' dv')} \lesssim \|F\|_{L^2}
\]
uniformly in $t \in (0,\pi/100]$.
Once this is shown, we see that $Z_{\bar{1}} \mathcal{S}f \in L^2(\Sigma_t)$ uniformly in $t \in (0,\pi/100]$, so $Z_{\bar{1}} \mathcal{S}f \in L^2(\Omega)$. Since $\mathcal{S}$ is real, we also have $Z_1 \mathcal{S}f \in L^2(\Omega)$. Thus $\mathcal{S}f \in S^{1,2}(\Omega)$, as desired.

To prove the claim, by Parseval on $L^2(\R^2 / \Lambda)$, we need only show that
\begin{equation} \label{eq:Kt_mulitplier_bdd}
\sup_{(m,n) \in \Z^2} \left|\iint_{\substack{|u| \leq \pi \\ |v| \leq \pi/2}} K_t(u,v) e^{-i(mu+nv)} du dv \right| \lesssim 1
\end{equation}
uniformly in $t \in (0,\pi/100]$. But we have $\Real K_t(-u,-v) = - \Real K_t(u,v)$ and $\Imag K_t(-u,v) = -\Imag K_t(u,v)$. One thus concludes the cancellation condition 
\[
\iint_{\norm{(u,v)} \leq r} K_t(u,v) du dv = 0
\]
for all $r$.
Next, for $|u| \leq \pi$ and $|v| \leq \pi/2$, a straightforward computation shows
\begin{align*}
 |K_t(u,v)| & \lesssim \norm{(u,v)}^{-3} , \\
 |\frac{\partial K_t}{\partial u}(u,v)| & \lesssim \norm{(u,v)}^{-4} , \\
 |\frac{\partial K_t}{\partial v}(u,v)| & \lesssim \norm{(u,v)}^{-5} ,
\end{align*}
all uniformly in $t \in (0,\pi/100]$. Combining these with the cancellation condition, one can obtain \eqref{eq:Kt_mulitplier_bdd} in a standard manner.
\end{proof}

To proceed further, define 
\[
\mathfrak{M}f(\zeta) := \sup_{t \in (0,\pi/100]} |\frac{d}{dt} \mathcal{S}f(\gamma_{\zeta}(t))|, \quad \zeta \in \Sigma.
\]
\begin{prop} \label{prop:frakM_bdd}
$\mathfrak{M}$ is bounded on $L^2(\Sigma,d\sigma)$.
\end{prop}

To prove this, write $M$ for the Hardy--Littlewood maximal operator adapted to the non-isotropic geometry on $\Sigma$:
\[
Mf(\zeta) := \sup_{t > 0} \frac{1}{|Q_t(\zeta)|} \int_{Q_t(\zeta)} |f(\eta)| d\sigma(\eta)
\]
where $Q_t(\zeta)$ is the non-isotropic ball given by  
\[
Q_t(\zeta) := \{\Phi(0,u,v) \in \Sigma \colon \norm{(u'-u,v'-v)} \leq t\} \quad \text{if $\zeta = \Phi(0,u',v')$.}
\]
(Recall $\norm{(u,v)} := \max\{|u|,|v|^{1/2}\}$.) It is known that $M$ is bounded on $L^2(\Sigma)$. Since we proved $\mathcal{N}$ is also bounded on $L^2(\Sigma)$, Proposition~\ref{prop:frakM_bdd} is a consequence of the following
\begin{prop} \label{prop:frakM_ptwise}
 For $f \in L^2(\Sigma)$, we have
\[
\mathfrak{M}f(\zeta) \leq C (M \mathcal{N}f(\zeta) + Mf(\zeta) )
\]
for $d\sigma$-almost every $\zeta \in \Sigma$.
\end{prop}


\begin{proof} [Proof of Proposition~\ref{prop:frakM_ptwise}]
In fact, if $\zeta = \Phi(0,u',v')$, then $\gamma_{\zeta}(t) = \Phi(t,u',v')$, so for $t > 0$ we can differentiate \eqref{eq:Sf_uv_nice} under the integral and get
\[
\frac{d}{dt} \mathcal{S}f(\gamma_{\zeta}(t)) = \frac{1}{8\pi} \iint_{\substack{|u| \leq \pi \\ |v| \leq \pi/2}} F(u'-u,v'-v) \frac{\partial k}{\partial t}(\frac{t}{\sqrt{2}},u,v) du dv ,
\]
where $F(u,v):= f(\Phi(0,u,v))$. We also have, from \eqref{eq:N_singular}, that
\[
\mathcal{N}f(\zeta)
= \lim_{\varepsilon \to 0^+} \frac{1}{8\pi} \iint_{\substack{|u| \leq \pi , |v| \leq \pi/2 \\ \norm{(u,v)} \geq \varepsilon}} F(u'-u,v'-v) \frac{\partial k}{\partial t}(0,u,v) du dv 
\]
and
\[
Mf(\zeta) \simeq \sup_{0 < t < 100} \frac{1}{t^3}  \iint_{\norm{(u,v)} \leq t} |F(u'-u,v'-v)| du dv.
\]
Using the bound
\[
|\frac{\partial k}{\partial t}(t,u,v)| \lesssim (|t|+|u|+|v|^{1/2})^{-3},
\]
it is easy to see that
\begin{equation}\label{max1}
\begin{split}
    \left| \frac{1}{8\pi} \iint_{\norm{(u,v)} \leq 10 t} F(u'-u,v'-v) \frac{\partial k}{\partial t} (\frac{t}{\sqrt{2}},u,v) du dv \right| \leq C Mf(\zeta);
\end{split}
\end{equation}
in fact, 
\begin{equation*}
 \left|\frac{1}{8\pi} \iint_{\norm{(u,v)} \leq 10 t} F(u'-u,v'-v) \frac{\partial k}{\partial t} (\frac{t}{\sqrt{2}},u,v) du dv \right| \leq
 C Mf(\zeta) \iint_{\norm{(u,v)} \leq 10 t} \frac{1}{
(|t|+|u|+|v|^{1/2})^{3}} du dv
\end{equation*}
and
 \[ \begin{split}
\iint_{\norm{(u,v)} \leq 10 t} \frac{1}{
(|t|+|u|+|v|^{1/2})^{3}} du dv
\leq 
\iint_{\norm{(u,v)} \leq 10 t} \frac{1}{
t^3} du dv
= C.
\end{split}
\]
Furthermore, using the bound
\[
|\frac{\partial k}{\partial t}(t,u,v) - \frac{\partial k}{\partial t}(0,u,v) | \lesssim \frac{t}{(|u|+|v|^{1/2})^4} \quad \text{when $\norm{(u,v)} \geq t$},
\]
one can show that
\begin{equation}\label{max3}
\begin{split}
    \left| \frac{1}{8\pi} \iint_{\substack{|u| \leq \pi \\ |v| \leq \pi/2 \\ \norm{(u,v)} > 10t}} F(u'-u,v'-v) \left( \frac{\partial k}{\partial t} (\frac{t}{\sqrt{2}},u,v) - \frac{\partial k}{\partial t} (0,u,v) \right) du dv  \right|
    \leq C M f(\zeta).
\end{split}
\end{equation}
This is because
\begin{multline*}
 \left|\frac{1}{8\pi} \iint_{\substack{|u| \leq \pi \\ |v| \leq \pi/2 \\ \norm{(u,v)} > 10t}} F(u'-u,v'-v) \left( \frac{\partial k}{\partial t} (\frac{t}{\sqrt{2}},u,v) - \frac{\partial k}{\partial t} (0,u,v) \right) du dv \right| \\
  \lesssim 
 Mf(\zeta) \iint_{\norm{(u,v)} > 10t} \frac{t}{(|u|+|v|^{1/2})^4} du dv \lesssim Mf(\zeta).
\end{multline*}
 
Finally, we claim that for almost every $\zeta \in \Sigma$,
\begin{equation}\label{max4}
\begin{split}
    \left| \frac{1}{8\pi}  \iint_{\substack{|u| \leq \pi \\ |v| \leq \pi/2 \\ \norm{(u,v)} > 10 t}} F(u'-u,v'-v) \frac{\partial k}{\partial t} (0,u,v) du dv \right|
    \leq & C (M\mathcal{N}f + Mf)(\zeta).
\end{split}
\end{equation}
In fact, the operators on both sides of the inequality are bounded on $L^2(\Sigma)$, and $C^{\infty}(\Sigma)$ is dense in $L^2(\Sigma)$, we may assume, without loss of generality, that $f \in C^{\infty}(\Sigma)$. Suppose now $\psi \in C^{\infty}_c(\Sigma)$ is supported in $Q_1(\Phi(0,0,0))$, with $\int_{\Sigma} \psi d\sigma = 1$. Let $\Psi_t(u,v) = \sqrt{2}^{-1} t^{-3} \psi(\Phi(0,t^{-1} u, t^{-2} v))$ for $t \in (0,\pi/100]$ and $(u,v) \in \R^2$, so that 
\[
\iint_{|u| \leq \pi, |v| \leq \pi/2} \Psi_t(u,v) du dv = 1.
\]
Recall $K(u,v) = p.v. \frac{\partial k}{\partial t}(0,u,v)$ and let 
\[
\tilde{K}_t(u,v) = K*\Psi_t(u,v)
\]
be the convolution of $K$ with $\Psi_t$ on $\R^2/\Lambda$. 
We will prove that
\begin{equation} \label{eq:partial_k_Ktilde_t}
|\frac{\partial k}{\partial t}(0,u,v) 1_{\norm{(u,v)} > 10 t} - \tilde{K}_t(u,v) | \lesssim \frac{t}{(|t|+|u|+|v|^{1/2})^4}.
\end{equation}
In fact, if  $\norm{(u,v)} > 10 t$, then the mean value theorem gives
\begin{align*}
| \frac{\partial k}{\partial t}(0,u,v) - \tilde{K}_t(u,v) |
&\leq \iint_{\norm{(u,v)} \leq t} \Big| \frac{\partial k}{\partial t}(0,u',v')- \frac{\partial k}{\partial t}(0,u'-u,v'-v)\Big| |\Psi_t(u,v)| du dv \\
&\lesssim \frac{t}{(|u|+|v|^{1/2})^4}.
\end{align*}
On the other hand, if $\norm{(u',v')} \leq 10 t$, then 
\begin{align*}
\tilde{K}_t(u',v') &= \iint_{\norm{(u,v)} \leq 100 t} \Big[ \Psi_t(u'-u,v'-v) - \Psi_t(u',v') \Big] \frac{\partial k}{\partial t}(0,u,v) du dv  \\
& \quad + \iint_{\substack{|u| \leq \pi, |v| \leq \pi/2 \\ \norm{(u,v)} > 100 t}} \Psi_t(u'-u,v'-v) \frac{\partial k}{\partial t}(0,u,v) du dv ,
\end{align*}
with the last integral being both zero because $\Phi_t(u'-u,v'-v) = 0$ there.
As a result, if $\norm{(u',v')} \leq 10 t$, then 
\begin{align*}
|\tilde{K}_t(u',v')| &\leq \iint_{\norm{(u,v)} \leq 100 t} \frac{1}{t^3}  \Big( \frac{|u|}{t} + \frac{|v|}{t^2} \Big) \frac{1}{\norm{(u,v)}^3} du dv \lesssim \frac{1}{t^3}.
\end{align*}
This proves \eqref{eq:partial_k_Ktilde_t}, and hence
\begin{multline*}
   \left| \frac{1}{8\pi}  \iint_{\substack{|u| \leq \pi \\ |v| \leq \pi/2 \\ \norm{(u,v)} > 10 t}} F(u'-u,v'-v) \frac{\partial k}{\partial t} (0,u,v) du dv \right| \\
   \leq  
      \left| \frac{1}{8\pi}  \iint_{\substack{|u| \leq \pi \\ |v| \leq \pi/2}} F(u'-u,v'-v) \tilde{K}_t(u,v) du dv \right| + Mf(\zeta).
\end{multline*}
But 
\[
\begin{split}
\frac{1}{8\pi} \iint_{\substack{|u| \leq \pi \\ |v| \leq \pi/2}} F(u'-u,v'-v) \tilde{K}_t(u,v) du dv &=  F*\frac{1}{8\pi} K*\Psi_t(u,v) \\
&= NF*\Psi_t(u,v),
\end{split}
\]
where $N F(u',v') = \mathcal{N}f(\Phi(0,u',v'))$.  Thus 
\[
\left| \frac{1}{8\pi}  \iint_{\substack{|u| \leq \pi \\ |v| \leq \pi/2}} F(u'-u,v'-v) \tilde{K}_t(u,v) du dv \right| \lesssim M \mathcal{N} f(\zeta).
\]
This proves \eqref{max4}.

\end{proof}

We are now ready to complete the proof of Theorem~\ref{thm:NeumannL2}. We need to show that if $f \in L^2(\Sigma)$ and $u = \mathcal{S}f$, then
\[
\begin{cases}
Lu = 0 \quad \text{on $\Omega$} \\
\nabla_{\nu} u = (-\frac{1}{2}I + \mathcal{N})f \quad \text{on  $\Sigma$}
\end{cases}
\]
in the weak sense, i.e.
\[
0 = \int_{\Omega} ( \langle \nabla_b u , \nabla_b \phi \rangle + \frac{R}{4} u \phi) + \int_{\Sigma} (-\frac{1}{2}I + \mathcal{N})f \phi
\]
for all $\phi \in C^{\infty}(\overline{\Omega})$. To do so, note that
\[
\int_{\Omega} (\langle \nabla_b u , \nabla_b \phi \rangle + \frac{R}{4} u \phi)
= \lim_{t \to 0^+} \int_{\Omega_t} (\langle \nabla_b u , \nabla_b \phi \rangle + \frac{R}{4} u \phi)
\]
and $u \in C^{\infty}(\overline{\Omega_t})$ for $t > 0$. Thus for $\phi \in C^{\infty}(\overline{\Omega})$ we can now integrate by parts, and using $Lu = 0$ inside $\Omega_t$, see that the above is equal to 
\[
-\lim_{t \to 0^+} \int_{\Sigma_t} \nabla_{\nu} u \, \phi \, d\sigma_t = -\lim_{t \to 0^+} \int_{\Sigma} \frac{d}{dt} u(\gamma_{\zeta}(t)) \phi(\gamma_{\zeta}(t)) \frac{d\sigma_t}{d\sigma} d\sigma.
\]
Our earlier estimate for $\mathfrak{M}$ allows one to show that
\begin{equation} \label{eq:ptwise}
\frac{d}{dt} u(\gamma_{\zeta}(t)) \to (-\frac{1}{2}I + \mathcal{N})f(\zeta) 
\end{equation}
for almost every $\zeta \in \Sigma$ as $t \to 0^+$: if $f \in L^2(\Sigma)$, then choose $g \in C^{\infty}(\Sigma)$ so that $\|f - g\|_{L^2(\Sigma)} < \varepsilon_0$. Then for any $\alpha > 0$,
\[
\begin{split}
&\{\zeta \in \Sigma \colon \limsup_{t \to 0^+} |\frac{d}{dt} u(\gamma_{\zeta}(t)) - (-\frac{1}{2}I + \mathcal{N})f(\zeta)| > \alpha\} \\
& \subseteq \{\zeta \in \Sigma \colon \limsup_{t \to 0^+} |\frac{d}{dt} \mathcal{S}(f-g)(\gamma_{\zeta}(t)) - (-\frac{1}{2}I + \mathcal{N})(f-g)(\zeta)| > \alpha\} \\
& \subseteq \{\zeta \in \Sigma \colon \mathfrak{M}(f-g)(\zeta) > \alpha/2 \} \cup \{ \zeta \in \Sigma \colon  |(-\frac{1}{2}I + \mathcal{N})(f-g)(\zeta)| > \alpha/2\}.
\end{split} 
\]
The measure of the last set with respect to $d\sigma$ is at most 
\[
C \alpha^{-2} \Big( \|\mathfrak{M}(f-g)\|_{L^2(\Sigma)}^2 +  \|(-\frac{1}{2} I + \mathcal{N}) (f-g)\|_{L^2(\Sigma)}^2 \Big) \lesssim \alpha^{-2} \|f-g\|_{L^2(\Sigma)}^2 \lesssim \alpha^{-2} \varepsilon_0^2.
\]
Letting $\varepsilon_0 \to 0^+$, we see that \eqref{eq:ptwise} holds $d\sigma$ a.e. on $\Sigma$.

Our earlier estimate for $\mathfrak{M}$ also allows us to further apply the dominated convergence theorem and show that 
\[
-\lim_{t \to 0^+} \int_{\Sigma} \frac{d}{dt} u(\gamma_{\zeta}(t)) \phi(\gamma_{\zeta}(t))  \frac{d\sigma_t}{d\sigma} d\sigma =
-\int_{\Sigma} (-\frac{1}{2}I + \mathcal{N})f \phi.
\]
This shows that $u$ solves the desired equation in the weak sense. 

For uniqueness, we need to show that if $u \in S^{1,2}(\Omega)$ and 
\[
0 = \int_{\Omega} \langle \nabla_b u , \nabla_b \phi \rangle + \frac{R}{4} u \phi \quad \text{for all $\phi \in C^{\infty}(\overline{\Omega})$},
\]
then $u = 0$ a.e. in $\Omega$. By density, the above condition actually implies 
\[
0 = \int_{\Omega} \langle \nabla_b u , \nabla_b \phi \rangle + \frac{R}{4} u \phi \quad \text{for all $\phi \in S^{1,2}(\Omega)$}.
\]
Setting $\phi = u$ and using $R > 0$, we obtain $u = 0$ a.e. on $\Omega$.

\section{Proof of Proposition~\ref{prop:clifford-steklov-regularity}} \label{sect6}

Let $\Omega$ be the domain in $\mathbb{S}^3$ defined by \eqref{eq:Omega_def}. Let $u \in S^{1,2}(\Omega)$ be a minimizer of $\mathcal{F}(u)$ under the condition that $\oint_M |u|^2 d\sigma = 1$. We want to show that $u \in C^\infty(\overline{\Omega})$. In fact this will be the case as long as $u$ is a critical point.

First, since
\[
0 = \frac{1}{2} \left. \frac{d}{dt} \right|_{t=0} \frac{\mathcal{F}(u+t\phi)}{\oint_M |\Tr(u+t\phi)|^2 d\sigma} = \int_{\Omega} \Big( \langle \nabla_b u , \nabla_b \phi \rangle + \frac{R}{4} u \phi \Big) - \mathcal{F}(u) \oint_{M} \Tr(u) \phi d\sigma
\]
for all $\phi \in C^{\infty}_c(\overline{\Omega})$, we have 
\begin{equation*}
 \begin{cases}
  Lu = 0 , & \text{in $\Omega$} , \\
  \nabla_\nu u = \mu \Tr(u) , & \text{on $\partial\Omega$} ,
 \end{cases}
\end{equation*}
where $\mu := -\mathcal{F}(u)$. By Theorem~\ref{thm:NeumannL2}, we have $u = \mathcal{S}f$ where $f := \mu (-\frac{1}{2} I + \mathcal{N})^{-1} \Tr(u)$. 

\begin{Lem} \label{lem:smoothing}
For any $s \geq 0$, we have
\[
\|\Tr \mathcal{S}f\|_{W^{s+\frac{1}{2},2}(\Sigma)} \leq C_s \|f\|_{W^{s,2}(\Sigma)} \quad \text{for all $f \in W^{s,2}(\Sigma)$}.
\]
\end{Lem}

Here we identify a function on $\Sigma$ with a function on the Lie group $\R^2 / \Lambda$, and $W^{s,2}(\Sigma)$ is the corresponding (isotropic) Sobolev space on $\R^2 / \Lambda$, so that 
\[
\|f\|_{W^{s,2}(\Sigma)} := \Big( \frac{1}{2\pi^2} \sum_{\substack{(m,n) \in \Z^2 \\ m \equiv n \pmod{2}}} (1+|m|^2+|n|^2)^s |\widehat{F}(m,n)|^2 \Big)^{1/2}
\]
if $F(u,v) := f(\Phi(0,u,v))$.

\begin{proof}[Proof of Lemma~\ref{lem:smoothing}]
To prove the Lemma, note that if $f \in L^2(\Sigma)$ and we write $F(u',v') = f(\Phi(0,u',v'))$ with $SF(u',v') := (\Tr \mathcal{S} f) (\Phi(0,u',v'))$, then
\[
SF = F*\frac{\sqrt{2}}{8\pi} K_0 \quad \text{a.e. on $\Sigma$}
\]
where $K_0(u',v') = k(0,u',v')$. The reason that this works for $f \in L^2(\Sigma)$ is because the identity holds for $f \in C^{\infty}(\Sigma)$, and that both sides are bounded linear maps on $L^2(\Sigma)$. Thus 
\[
\widehat{SF}(m,n) = \widehat{F}(m,n) \frac{\sqrt{2}}{8\pi} \widehat{K_0}(m,n)
\]
for all $(m,n) \in (\frac{1}{2\pi} \Lambda)^*$.
It is easy to see that 
\[
|\widehat{K_0}(m,n)| \lesssim (1+|m|+|n|^{1/2})^{-1}.
\]
In fact, since $K_0 \in L^1(\Sigma)$, we have $|\widehat{K_0}(m,n)| \lesssim 1$. Also, since 
\[
|K_0(u,v)| \lesssim \norm{(u,v)}^{-2}, \quad |\frac{\partial^2 K_0}{\partial u^2}| + |\frac{\partial K_0}{\partial v}| \lesssim \norm{(u,v)}^{-4},
\]
by writing
\[
\widehat{K_0}(m,n) = \iint_{\norm{(u,v)} \leq (|m|+|n|^{1/2})^{-1}} +  \iint_{\substack{|u| \leq \pi, |v| \leq \pi/2 \\ \norm{(u,v)} \geq (|m|+|n|^{1/2})^{-1}}} K_0(u,v) e^{-i(mu+nv)} du dv,
\]
and integrating by parts in the second integral using either $e^{-imu} = -\frac{1}{m^2} \frac{\partial^2}{\partial u^2} e^{-imu}$ or $e^{-inv} = \frac{i}{n} \frac{\partial}{\partial v} e^{-inv}$, we obtain $|\widehat{K_0}(m,n)| \lesssim (|m|+|n|^{1/2})^{-1}$. This proves our bound for $\widehat{K_0}(m,n)$; in fact, we will only use the weaker bound 
\[
|\widehat{K_0}(m,n)| \lesssim (1+|m|+|n|)^{-1/2}.
\]
Now if $f \in W^{s,2}(\Sigma)$, then this weaker estimate for $\widehat{K}(m,n)$ allows us to bound
\begin{align*}
\|\Tr \mathcal{S} f\|_{W^{s+1/2,2}(\Sigma)} &
\simeq \Big( \sum_{\substack{(m,n) \in \Z^2 \\ m \equiv n \pmod{2}}} (1+|m|^2+|n|^2)^{s+\frac{1}{2}} |\widehat{SF}(m,n)|^2 \Big)^{1/2} \\
&\simeq \Big( \sum_{\substack{(m,n) \in \Z^2 \\ m \equiv n \pmod{2}}} (1+|m|^2+|n|^2)^{s+\frac{1}{2}} |\widehat{F}(m,n) \widehat{K_0}(m,n)|^2 \Big)^{1/2} \\
&\lesssim \Big( \sum_{\substack{(m,n) \in \Z^2 \\ m \equiv n \pmod{2}}} (1+|m|^2+|n|^2)^s |\widehat{F}(m,n)|^2 \Big)^{1/2} \\
&\simeq \|f\|_{W^{s,2}(\Sigma)}. \qedhere
\end{align*}
\end{proof}

Recall now our critical point $u = \mathcal{S}f$ where $f := \mu (-\frac{1}{2}I + \mathcal{N})^{-1} \Tr(u)$. By Lemma~\ref{lem:smoothing}, for all $s \geq 0$, we have
\[
\|\Tr u\|_{W^{s+\frac{1}{2},2}(\Sigma)} \lesssim \|f\|_{W^{s,2}(\Sigma)}.
\]
Using the bound from Section~\ref{sect4}, we have
\[
\|f\|_{W^{s,2}(\Sigma)} \lesssim \mu \|\Tr(u)\|_{W^{s,2}(\Sigma)}.
\]
Putting these together, we see that 
\[
\|\Tr(u)\|_{W^{k,2}(\Sigma)} \lesssim \mu^{2k} \|\Tr(u)\|_{L^2(\Sigma)} \lesssim \mu^{2k} \|u\|_{S^{1,2}(\Omega)}
\]
for all $k \in \mathbb{N}$. Thus $\Tr(u) \in C^{\infty}(\Sigma)$, and hence $\nabla_{\nu} u \in C^{\infty}(\Sigma)$. As a result, $f = \mu (-\frac{1}{2}I + \mathcal{N})^{-1} \nabla_{\nu} u \in C^{\infty}(\Sigma)$. Thus $u = \mathcal{S} f \in C^{\infty}(\overline{\Omega})$, as desired.

\section{The boundary Yamabe contact form on $\Sigma$}
\label{sec:compute-boundary-form}

We conclude this article by proving our results about the boundary Yamabe constant.

First we prove that the boundary Yamabe constant $Y(M,T^{1,0})$ is finite if and only if the first Dirichlet eigenvalue of the CR Yamabe operator is positive.

\begin{prop}
 \label{finite-yamabe-to-dirichlet-eigenvalue}
 Let $(M^3,T^{1,0},\theta)$ be a closed pseudohermitian manifold with boundary having no characteristic points.
 Then $Y(M,T^{1,0}) > -\infty$ if and only if $\lambda_{1,D}(L) > 0$.
\end{prop}

\begin{proof}
 Set $\lambda := \lambda_{1,D}(L)$ and let $u$ be the first Dirichlet eigenvalue;
 i.e.\ $u$ is the unique nonnegative function in $S^{1,2}(M)$ such that $\int_M u^2\,\theta \wedge d\theta = 1$ and
 \begin{equation*}
  \begin{cases}
   Lu = \lambda u , & \text{in $M$}, \\
   u = 0 , & \text{on $\partial M$} .
  \end{cases}
 \end{equation*}
 Jerison proved~\cite[Theorem~7.1]{Jerison1981} that $u \in C^\infty(M)$.
 Bony's maximum principle~\cite[Corollaire~3.1]{Bony1969} then implies that $u>0$.
 Note that $\mc{F}(u) = \lambda$.
 Additionally, for any fixed $v \in C^\infty(M)$, it holds that
 \begin{equation}
  \label{eqn:perturb}
  \begin{split}
   \mc{F}(v + tu) & = \lambda t^2 + \lambda t \int_M uv \, \theta \wedge d\theta + t \oint_{\partial M} v \, Bu \, d\sigma + \mc{F}(v) , \\
   \oint_{\partial M} \lvert v+tu \rvert^3 \, d\sigma & = \oint_{\partial M} \lvert v \rvert^3 \, d\sigma .
  \end{split}
 \end{equation}
 
 Suppose first that $\lambda < 0$.
 For any fixed $v \in C^\infty(M)$, we deduce from~\eqref{eqn:perturb} that $\mc{F}(v+tu) \to -\infty$ as $t \to \infty$.
 Therefore $Y(M,T^{1,0}) = -\infty$.
 
 Suppose next that $\lambda = 0$.
 Fix a $v \in C^\infty(M)$ such that $v\rvert_{\partial M} = -Bu$.
 We deduce from~\eqref{eqn:perturb} that $\mc{F}(v+tu) \to -\infty$ as $t \to \infty$.
 Therefore $Y(M,T^{1,0}) = -\infty$.
 
 Suppose finally that $\lambda>0$.
 Fix a constant $c \in (0,\lambda)$.
 A standard variational argument implies the existence of a weak solution $u \in S^{1,2}(M)$ of
 \begin{equation*}
  \begin{cases}
   Lu = cu , & \text{in $M$}, \\
   u = 1 , & \text{on $\partial M$}.
  \end{cases}
 \end{equation*}
 Combining Jerison's regularity results~\cite[Theorem~7.1]{Jerison1981} with Bony's maximum principle~\cite[Corollaire~3.1]{Bony1969} implies that $u \in C^\infty(M)$ is positive.
 Set $\widehat{\theta} := u^2\theta$.
 By the conformal covariance of the CR Yamabe operator,
 \begin{equation*}
  R^{\widehat{\theta}} = 4L^{u^2\theta}(1) = 4cu^{-2} > 0 .
 \end{equation*}
 In particular, there are constants $C_1,C_2>0$ such that
 \begin{equation*}
  \mathcal{F}^{\widehat{\theta}}(v) \geq \int_M \left( \lvert \nablab v \rvert^2 + C_1 \lvert v \rvert^2 \right) \, \widehat{\theta} \wedge d\widehat{\theta} - C \oint_{\partial M} \lvert v\rvert^2 \, d\sigma_{\widehat{\theta}}
 \end{equation*}
 for all $v \in C^\infty(M)$.
 On the one hand, H\"older's inequality implies that $\oint \lvert v \rvert^2 \lesssim \bigl( \oint \lvert v \rvert^3 \bigr)^{2/3}$.
 On the other hand, the Sobolev trace inequality~\cite[Theorem~1.4]{Nhieu2001} and a standard partition of unity argument implies that
 \begin{equation*}
  \int_M \left( \lvert \nablab v \rvert^2 + C_1 \lvert v \rvert^2 \right) \, \widehat{\theta} \wedge d\widehat{\theta} \gtrsim \left( \oint_{\partial M} \lvert v \rvert^3 \, d\widehat{\sigma} \right)^{\frac{2}{3}} .
 \end{equation*}
 Therefore $Y(M,T^{1,0}) > -\infty$.
%
\end{proof}

Second we prove that $\mu_{1}(L)$ and $Y(M,T^{1,0})$, if finite, have the same sign provided that minimizers of the former are smooth.

\begin{prop}
 \label{prop:same-signs}
  Let $(M^3,T^{1,0},\theta)$ be a closed pseudohermitian manifold with boundary having no characteristic points.
  Suppose additionally that the minimizer of~\eqref{eqn:first-neumann} is smooth. 
 Then $\mu_{1}(L)$ and $Y(M,T^{1,0})$ have the same sign.
\end{prop}

\begin{proof}
 Denote $Y := Y(M,T^{1,0})$ and $\mu_{1} := \mu_{1}(L)$.
 
 Suppose first that $Y > 0$ (resp.\ $Y \geq 0$).
 H\"older's inequality implies that if $u \in C^\infty(M)$, then
 \begin{equation*}
  \mc{F}(u) \geq Y \left( \oint_{\partial M} \lvert u\rvert^3 \, d\sigma \right)^{\frac{2}{3}} \geq Y \mathrm{Vol}(\partial M)^{-\frac{1}{3}} \oint_{\partial M} \lvert u \rvert^2 \, d\sigma .
 \end{equation*}
 Therefore $\mu_1 > 0$ (resp.\ $\mu_1 \geq 0$).
 
 Suppose next that $\mu_1 > 0$ (resp.\ $\mu_{1} \geq 0$).
 By restricting to nonnegative functions and applying the Sobolev trace embedding theorem~\cite[Theorem~1.4]{Nhieu2001}, we deduce that there is a nonnegative minimizer $u \in S^{1,2}(M)$ of $\mu_{1}$.
 By assumption, $u$ is smooth.
 A straightforward computation implies that $u$ is a weak solution of
 \begin{equation}
  \label{eqn:steklov}
  \begin{cases}
   Lu = 0 , & \text{in $M$} , \\
   Bu = \mu_{1}u , & \text{on $\partial M$} .
  \end{cases}
 \end{equation}
 The strong maximum principle~\cite[Corollaire~3.1]{Bony1969} implies that $u$ is positive in the interior of $\Omega$.
 The Hopf Lemma~\cite[Corollary~2.1]{Monticelli2010} implies that $u$ is also positive on the boundary.
 Therefore $u>0$.
 Set $\widehat{\theta} := u^2\theta$.
 Equation~\eqref{eqn:steklov} implies that $H^{\widehat{\theta}} > 0$ (resp.\ $H^{\widehat{\theta}} \geq 0$).
 By conformal covariance,
 \begin{equation*}
  \mc{F}^\theta(u^{-1}v) = \mc{F}^{\widehat{\theta}}(v) = 2\int_M \lvert \nabla_b v \rvert^2 \, \widehat{\theta} \wedge d\widehat{\theta} + \frac{1}{3} \oint_M H^{\widehat{\theta}} \lvert v \rvert^2 \, \widehat{d\sigma} .
 \end{equation*}
 We conclude from the Sobolev trace embedding theorem that $Y>0$ (resp.\ $Y \geq 0$).
%
\end{proof}

Finally, we construct a scalar flat contact form on $\Omega$ with respect to which $\Sigma$ has constant $p$-mean curvature.

\begin{proof}[Proof of Theorem~\ref{thm:clifford-extremals}]
 Introduce new coordinates $(r,u,v) \in [0,\pi] \times [0,2\pi] \times [0,2\pi]$ on $\Omega$ by
 \begin{align*}
  \zeta^1 & = e^{i(u+v)}\sin(r/2) , \\
  \zeta^2 & = e^{i(u-v)}\cos(r/2) ,
 \end{align*}
 making the obvious identifications.
 A straightforward computation yields
 \begin{equation*}
  Z_1 = e^{-2iu}\left( \partial_r + \frac{1}{2i}\cos(r)\partial_u + \frac{1}{2i} \csc(r)\partial_v \right) .
 \end{equation*}
 It readily follows that
 \begin{equation*}
  L = -2\partial_r^2 - \frac{1}{2}\cot^2(r)\partial_u^2 - \frac{1}{2}\csc^2(r)\partial_v^2 - \cot(r)\csc(r)\partial_{uv}^2 - 2\cot(r)\partial_r + \frac{1}{2} .
 \end{equation*}
 We conclude that a smooth function $u=u(r)$ on $\Omega$ solves $Lu = 0$ if and only if
 \begin{equation}
  \label{eqn:clifford-ode}
  \begin{cases}
  \partial_r^2 u + \cot(r) \partial_r u - \frac{1}{4}u = 0 , \\
  \displaystyle\lim_{r\to\pi} u(r) < \infty .
  \end{cases}
 \end{equation}
 By the change of variables $x=\cos r$, we see that \eqref{eqn:clifford-ode} is equivalent to
 \begin{equation}
  \label{eqn:clifford-ode-x}
  \begin{cases}
  (1-x^2)\partial_x^2 u - 2x\partial_x u - \frac{1}{4}u = 0 , \\
  \displaystyle\lim_{x\to-1} u(x) < \infty .
  \end{cases}
 \end{equation}
 It is well-known~\cite[\S14.2 and \S14.8]{NIST:DLMF} that the space of solutions of \eqref{eqn:clifford-ode-x} is spanned by the Legendre function
 \begin{equation*}
  P_{-\frac{1}{2}}(x) = {}_2F_1 \left( \frac{1}{2} , \frac{1}{2} ; 1 ; \frac{1+x}{2} \right) .
 \end{equation*}
 Rewriting this in terms of $\lvert z^2\rvert^2 = \cos^2\frac{r}{2}$ and normalizing so that $u=1$ along $\Sigma$ yields the claimed solution of $Lu=0$.
 In particular, $u^2\theta$ is Webster-flat.
 Since $u$ depends only on $r$, we see that both $u$ and $Bu = \partial_ru$ are constant on $\Sigma$.
 Therefore $\Sigma$ has constant $p$-mean curvature with respect to $u^2\theta$.
\end{proof}

\section{Appendix}
\subsection{Appendix A}
The Mathematica codes used in Corollaries \ref{cor:khat1} and \ref{cor:khat2} are reproduced below for reference. In order to simplify expressions and outputs in Mathematica, some definitions in the codes are different from the definitions taken in the rest of the paper. We first point out the relations between these two sets of definitions and then list the codes:
\begin{align*}
&\e_0(a,u,v)=\frac{e_0(a)}{e_1(a)^{5/2}}\mathtt{error0[u, v]},\qquad \e_u(a,u,v)=\frac{e_0(a)}{e_1(a)^{7/2}}\mathtt{erroru1[u, v]} +\frac{1}{e_1(a)^{5/2}} 
\mathtt{erroru2[u, v]},
\\
&\e_v(a,u,v)=\frac{e_0(a)}{e_1(a)^{7/2}}\mathtt{errorv1[u, v]}+\frac{1}{e_1(a)^{5/2}}\mathtt{errorv2[u, v]},
\\
&\e_{uu}(a,u,v)=\frac{e_0(a)}{e_1(a)^{9/2}}\mathtt{erroruu1[u, v]}+\frac{1}{e_1(a)^{7/2}}\mathtt{erroruu2[u, v]},
\\
&\eps_{I0}=\frac{e_0(2)}{e_1(2)^{5/2}}\mathtt{epsiloni[0]},\qquad\eps_{I1}=\frac{e_0(2)}{e_1(2)^{5/2}}\mathtt{epsiloni[1]},
\\
&\eps_{I2}=\frac{e_0(2)}{e_1(2)^{7/2}}\mathtt{epsiloni[2,1]}+\frac{1}{e_1(2)^{5/2}}\mathtt{epsiloni[2,2]},
\\
&\eps_{I3}=\frac{e_0(2)}{e_1(2)^{9/2}}\mathtt{epsiloni[3,1]}+\frac{1}{e_1(2)^{7/2}}\mathtt{epsiloni[3,2]},
\\
&\eps_{J0}=\frac{e_0(1)}{e_1(1)^{5/2}}\mathtt{epsilonj[0]},\qquad\eps_{J1}=\frac{e_0(1)}{e_1(1)^{5/2}}\mathtt{epsilonj[1]},
\\
&\eps_{J2}=\frac{e_0(1)}{e_1(1)^{7/2}}\mathtt{epsilonj[2,1]}+\frac{1}{e_1(1)^{5/2}}\mathtt{epsilonj[2,2]}.
\end{align*}

\begin{verbatim}
d[u_, v_] = ((Cos[u] - Cos[v])^2 + Sin[v]^2)^(1/4);
n[u_, v_] = Sin[u]*Sin[v];
K[u_, v_] = n[u, v]/d[u, v]^6;
nu[u_, v_] = D[K[u, v], u]*d[u, v]^(10);
nv[u_, v_] = D[K[u, v], v]*d[u, v]^(10);
nuu[u_, v_] = D[K[u, v], u, u]*d[u, v]^(14);
d0[u_, v_] = (u^4/4 + v^2)^(1/4);
n0[u_, v_] = u*v;
K0[u_, v_] = n0[u, v]/(d0[u, v])^6;
n0u[u_, v_] = v (-5/4 u^4 + v^2);
K0u[u_, v_] = n0u[u, v]/d0[u, v]^(10);
n0v[u_, v_] = u (u^4/4 - 2 v^2);
K0v[u_, v_] = n0v[u, v]/d0[u, v]^(10);
n0uu[u_, v_] = u^3 (-15/2*v^3 + 15/8*u^4*v);
K0uu[u_, v_] = n0uu[u, v]/d0[u, v]^(14);
error0[u_, v_] := (6/4)*Abs[n0[u, v]]/d0[u, v]^4;
erroru1[u_, v_] := (10/4)*Abs[n0u[u, v]]/d0[u, v]^8;
erroru2[u_, v_] := (1451/3072*Abs[u^6*v] + Abs[u^2 v^3]/2 + Abs[v^4]/3)/d0[u, v]^(10);
errorv1[u_, v_] := (10/4)*Abs[n0v[u, v]]/d0[u, v]^8;
errorv2[u_, v_] := (683/15360 Abs[u^7] + Abs[u^3 v^2] + 7/6 Abs[u v^3])/d0[u, v]^(10);
erroruu1[u_, v_] := (14/4)*Abs[n0uu[u, v]]/d0[u, v]^(12);
erroruu2[u_, v_] := (1195/2048 Abs[u^9 v] + Abs[u^5 v^3]/8 + 5 Abs[u^3 v^4]/2 + 
     15 Abs[u v^5]/4)/d0[u, v]^(14);
norm[a_, u_, v_] := Max[Abs[u]/a, Sqrt[Abs[v]]];
region0[a_] := ImplicitRegion[norm[a, u, v] < 1, {u, v}];
region1[a_] := ImplicitRegion[norm[a, u, v] > 1, {u, v}];
b = 2/5; a1 = 2;
i[0] = FullSimplify[Integrate[Abs[u*v*K0[u, v]], {u, v} \[Element] region0[a1]]] 
i[1] = FullSimplify[2*Integrate[Abs[K0[a1, v]], {v, -1, 1}]]
i[2] = FullSimplify[2*Integrate[Abs[K0u[a1, v]], {v, -1, 1}]]
i[3] = FullSimplify[Integrate[Abs[K0uu[u, v]], {u, v} \[Element] region1[a1]]]
region0b[a1] = ImplicitRegion[norm[a1, u, v] < b, {u, v}];
e0[a1] = ToRadicals[MaxValue[(u^6/24 + u^2 v^2/2 + 2 u^6*(a1*b)^2/Factorial[8] + 
      Abs[v]^3*b^2/12)/d0[u, v]^(6), {u, v} \[Element] region0b[a1]]]
d0[a1*b, b^2]^2
epsiloni[0] = FullSimplify[Integrate[Abs[u*v]*error0[u, v], {u, v} \[Element] region0[a1]]] 
epsiloni[1] = FullSimplify[2*Integrate[error0[a1, v], {v, -1, 1}]]
epsiloni[2, 1] = FullSimplify[2*Integrate[erroru1[a1, v], {v, -1, 1}]] 
epsiloni[2, 2] = FullSimplify[2*Integrate[erroru2[a1, v], {v, -1, 1}]] 
epsiloni[3, 1] = FullSimplify[2*a1*Integrate[erroruu1[a1, v], {v, -1, 1}] +   
    4*Integrate[erroruu1[u, 1], {u, -a1, a1}]]
epsiloni[3, 2] = FullSimplify[2*a1*Integrate[erroruu2[a1, v], {v, -1, 1}] + 
    4*Integrate[erroruu2[u, 1], {u, -a1, a1}]]
epsiloni[4] = SetAccuracy[
 2*NIntegrate[Abs[nuu[u, v]]/d[u, v]^(14), {u, -a1*b, a1*b}, {v, b^2, Pi/2}, 
    Method -> "LocalAdaptive", AccuracyGoal -> 8, 
    IntegrationMonitor :> ((errors1 = Through[#1@"Error"]) &)]
+2*NIntegrate[Abs[nuu[u, v]]/d[u, v]^(14), {u, a1*b, Pi}, {v, -Pi/2, Pi/2}, 
    Method -> "LocalAdaptive", AccuracyGoal -> 8, 
    IntegrationMonitor :> ((errors2 = Through[#1@"Error"]) &)], 10]
SetAccuracy[Total@errors1 + Total@errors2, 10]
a2 = 1; 
j[0] = FullSimplify[Integrate[Abs[u*v*K0[u, v]], {u, v} \[Element] region0[a2]]]
j[1] = 2*Integrate[Abs[K0[u, 1]], {u, -a2, a2}]
j[2] = Integrate[Abs[K0v[u, v]], {u, v} \[Element] region1[a2]]
region0b[a2] = ImplicitRegion[norm[a2, u, v] < b, {u, v}];
e0[a2] = ToRadicals[MaxValue[(u^6/24 + u^2 v^2/2 + 2 u^6*(a2*b)^2/Factorial[8] + 
      Abs[v]^3*b^2/12)/d0[u, v]^(6), {u, v} \[Element] region0b[a2]]]
d0[a2*b, b^2]^2
epsilonj[0] = FullSimplify[Integrate[Abs[u*v]*error0[u, v], {u, v} \[Element] region0[a2]]]
epsilonj[1] = FullSimplify[2*Integrate[error0[u, 1], {u, -a2, a2}]]
epsilonj[2, 1] = FullSimplify[2*a2*Integrate[errorv1[a2, v], {v, -1, 1}] + 
   4*Integrate[errorv1[u, 1], {u, -a2, a2}]]
epsilonj[2, 2] = FullSimplify[2*a2*Integrate[errorv2[a2, v], {v, -1, 1}] + 
   4*Integrate[errorv2[u, 1], {u, -a2, a2}]]
SetAccuracy[
 2*NIntegrate[Abs[nv[u, v]]/d[u, v]^(10), {u, -a2*b, a2*b}, {v, b^2, Pi/2}, 
    Method -> "LocalAdaptive", AccuracyGoal -> 8, 
    IntegrationMonitor :> ((errors1 = Through[#1@"Error"]) &)] 
+2*NIntegrate[Abs[nv[u, v]]/d[u, v]^(10), {u, a2*b, Pi}, {v, -Pi/2, Pi/2}, 
    Method -> "LocalAdaptive", AccuracyGoal -> 8, 
    IntegrationMonitor :> ((errors2 = Through[#1@"Error"]) &)], 10]
SetAccuracy[Total@errors1 + Total@errors2, 10]
\end{verbatim}

\subsection{Appendix B}
We use the following codes in the Mathematica to verify the inequality
$|\widehat{K}(m,n)| < 4\pi$ when $1 \leq m \leq 3$, $1 \leq n \leq 7$, $m, n\in \mathbb{Z}$:
\begin{verbatim}
m := 1; For[n = 1, n <= 7, n += 2, 
 Print[NIntegrate[K[u, v] Sin[m*u] Sin[n*v], {u, -Pi, Pi}, {v, -Pi/2, Pi/2}, 
   Method -> "LocalAdaptive", AccuracyGoal -> 8, 
   IntegrationMonitor :> ((errors = Through[#1@"Error"]) &)]]; 
 Print[Total@errors]]
 
m := 2; For[n = 2, n <= 6, n += 2, 
 Print[NIntegrate[K[u, v] Sin[m*u] Sin[n*v], {u, -Pi, Pi}, {v, -Pi/2, Pi/2}, 
   Method -> "LocalAdaptive", AccuracyGoal -> 8, 
   IntegrationMonitor :> ((errors = Through[#1@"Error"]) &)]]; 
 Print[Total@errors]]
m := 3; For[n = 1, n <= 7, n += 2, 
 Print[NIntegrate[K[u, v] Sin[m*u] Sin[n*v], {u, -Pi, Pi}, {v, -Pi/2, Pi/2}, 
   Method -> "LocalAdaptive", AccuracyGoal -> 8, 
   IntegrationMonitor :> ((errors = Through[#1@"Error"]) &)]]; 
 Print[Total@errors]]
\end{verbatim}

\bibliographystyle{alpha}
\bibliography{References}
\end{document}